\documentclass[msom,nonblindrev]{informs3-modified}

\OneAndAHalfSpacedXI

\usepackage{natbib}
\usepackage{placeins}
\usepackage{cleveref}
\usepackage{subfigure}
\usepackage{enumerate}
\usepackage{accents}
\usepackage{amsmath}
\usepackage{booktabs}
\usepackage{booktabs}
\usepackage{float}
\usepackage{color}
\usepackage{xcolor}
\usepackage{soul}
\usepackage{algorithm}
\usepackage{algpseudocode}
\usepackage{mathtools}
\usepackage{bm}

\newcommand{\xx}{\boldsymbol{\mathit{x}}}

\newcolumntype{C}[1]{>{\centering\arraybackslash}p{#1}}

\newcommand{\vwp}{y}

\newcommand{\scriptG}{\mathcal{G}}

\crefformat{section}{\S#2#1#3}
\crefformat{subsection}{\S#2#1#3}
\crefformat{subsubsection}{\S#2#1#3}

 \bibpunct[, ]{(}{)}{,}{a}{}{,}%
 \def\bibfont{\small}%

\definecolor{var_color}{rgb}{0,0,0}

\TheoremsNumberedThrough     

\EquationsNumberedThrough    

\MANUSCRIPTNO{(MSOM-22-558-R1)} 

\begin{document}


\RUNAUTHOR{Hassanzadeh, Hosseini, Turner}

\RUNTITLE{Beyond Suspension: A Two-phase Methodology for Concluding Sports Leagues}

\TITLE{Beyond Suspension: A Two-phase Methodology for Concluding Sports Leagues}

\ARTICLEAUTHORS{%
\AUTHOR{Ali Hassanzadeh$^\text{a}$, Mojtaba Hosseini$^\text{b}$, John G. Turner$^\text{c}$}
\AFF{$^\text{a}$Alliance Manchester Business School, The University of Manchester, UK; $^\text{b}$ Tippie College of Business, University of Iowa, Iowa, USA; $^\text{c}$The Paul Merage School of Business, University of California, Irvine, California, USA }
\AFF{\EMAIL{ali.h@manchester.ac.uk}, \EMAIL{mojtaba-hosseini@uiowa.edu}, \EMAIL{john.turner@uci.edu}
}
} 

\ABSTRACT{%
\textbf{Problem definition:} Professional sports leagues may be suspended due to various reasons such as the recent COVID-19 pandemic. A critical question the league must address when re-opening is how to appropriately select a subset of the remaining games to conclude the season in a shortened time frame. 
\textbf{Academic/practical relevance:} Despite the rich literature on scheduling an entire season starting from a blank slate, concluding an existing season is quite different. Our approach attempts to achieve team rankings similar to that which would have resulted had the season been played out in full.
\textbf{Methodology:} We propose a data-driven model which exploits predictive and prescriptive analytics to produce a schedule for the remainder of the season comprised of a subset of originally-scheduled games.
Our model introduces novel rankings-based objectives within a stochastic optimization model, whose parameters are first estimated using a predictive model.
We introduce a deterministic equivalent reformulation along with a tailored Frank-Wolfe algorithm to efficiently solve our problem, as well as a robust counterpart based on min-max regret.
\textbf{Results:}
We present simulation-based numerical experiments from previous National Basketball Association (NBA) seasons 2004--2019, and show that our models are computationally efficient, outperform a greedy benchmark that approximates a non-rankings-based scheduling policy, and produce interpretable results.
\textbf{Managerial implications:} Our data-driven decision-making framework may be used to produce a shortened season with 25-50\% fewer games while still producing an end-of-season ranking similar to that of the full season, had it been played.
}



\maketitle

%



\section {Introduction}
\label{sec:Intro}

The novel coronavirus SARS-CoV-2, the causative agent of \mbox{COVID-19}, emerged in December 2019, causing the World Health Organization (WHO) to declare a pandemic.
In the ensuing months, governments rushed to enforce unprecedented quarantines, lockdowns, travel restrictions, and social distancing measures to reduce viral transmission
\citep{kaplan2020om}. As a result, supply chains were disrupted \citep{betcheva2020supply}, as well as many social activities.
Worldwide, professional sports leagues ceased activity, with the National Basketball Association (NBA) being the first in the US to suspend games, pausing the season as of March 11, 2020 \citep{nba-covid}. As the league  contemplated re-opening, a critical question was how to 
conclude the season in a shortened time frame.

Professional sports leagues globally adopt various formats including round-robin and elimination tournaments, with leagues such as soccer solely using round-robin, while major sports leagues in North America including the NBA, National Football League (NFL), National Hockey League (NHL), and Major League Baseball (MLB) combine a round-robin season with an elimination postseason (``playoffs").
Schedule adjustments are often required when a sports league is suspended (i.e, paused) for any reason (e.g., player strikes). The suspended 2019-20 season due to the \mbox{COVID-19} pandemic is a recent example.
Table \ref{table:suspensions} lists the number of suspensions in the NBA, NFL, NHL, and MLB from 1946--2021, along with the number of suspensions that led to shortened seasons. In summary, suspensions are an occasional occurrence of significant consequence, and our methodology can be applied to resuming a league's regular season from any such suspension.
\begin{table}[ht]
\centering
\footnotesize
\begin{tabular}{ccccc}
\toprule
League           & NBA & NFL & NHL & MLB \\\midrule
Number of Suspensions       & 6   & 8   & 4   & 9   \\
Suspensions with Shortened Seasons & 4   & 2   & 2   & 3\\\bottomrule
\end{tabular}
\vspace{0.1cm}
\caption{The number of disruptions to regular season games.}
\label{table:suspensions}
\end{table}
\vspace{-0.5cm}

After the NBA's \mbox{COVID-19} suspension in March 2020, there was much media speculation on the possible  actions the NBA would take: whether the 2019--20 regular season would be resumed, and how it would be concluded. 
The main directions the NBA could have taken were as follows:
\begin{enumerate}
    \item \textbf{Cancel Everything.} Cancel the remaining regular-season games as well as the playoffs. Determine the champion by vote.
    \item \textbf{Skip to Playoffs.} Cancel the remaining regular-season games. The top-ranked teams as of the suspension date qualify for the playoffs, which begin immediately.
    \item \textbf{Full Season.} Resume the regular season and play all 259 scheduled games, followed by the playoffs as usual. The season will end late in the year.
    \item \textbf{Shortened Season.} Select a subset of the remaining 259 regular-season games to be played, followed by the playoffs.
\end{enumerate}
\begin{figure}[ht]
    \centering
    \includegraphics[width=0.6\textwidth]{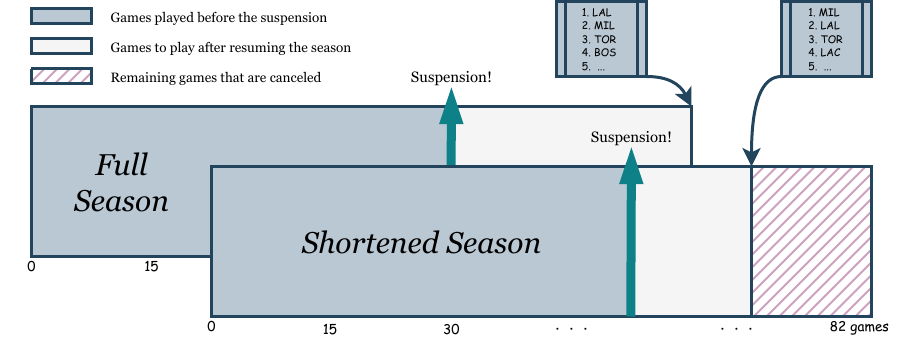}
    \caption{Two strategies to conclude the league: full season vs. shortened season after resuming the league}
    \label{fig:two-strategies}
\end{figure}
Options 1 and 2 would be unfair for several reasons. Indeed, teams may differ in the number of games played by the suspension date, along with the relative difficulties of the opponents they played.
Options 3 and 4, which we compare extensively throughout the paper, are illustrated in Figure \ref{fig:two-strategies}. Note that the rankings of the teams (by number of wins over the played games) depends on the specific set of games that are played, as well as the outcomes of those games.

Our focus in this paper is to propose a method which chooses a subset of games to conclude a shortened season that remains an asymmetrical round-robin tournament while producing end-of-season rankings that are as close as possible to the rankings that would result had the full season been played (i.e., no games canceled). There are, of course, many considerations that come into play when constructing a sports schedule. While our approach is less detailed than constructing a full timetable that incorporates not only the set of games to play but also their sequence (and corresponding travel schedule), our model is sufficiently general to allow for logical constraints on the subset of games chosen, which may be driven by specific practical considerations. 

At a high level, our model selects which games to include in the remainder of the season. To make these decisions, we develop several model components. First, we use a predictive model to estimate the likelihood of each game outcome for all games in the season that have not yet been played. Then, we generate one or more scenarios, and for each scenario we produce a ``target ranking" which ranks the relative performance of the teams if all games in the remainder of the season were played. We then use a prescriptive model to select a subset of games so that, in expectation, the ranking we get from our shortened season is as close to the target ranking as possible. 

Our paper is organized as follows. We review related works in \cref{sec:lit-review}, and mathematically define our problem in \cref{sec:problem_background}. In \cref{sec:models} we introduce our predictive and prescriptive models (binary classifiers and stochastic optimization models, respectively), and in \cref{sec:solution-methods} we develop solution techniques (one based on the Frank-Wolfe algorithm) for solving our prescriptive models efficiently. In \cref{sec:comp-exp} we use Monte Carlo simulation to  evaluate our models, and show that not only do the shortened seasons  produced by our best model have rankings that are close to the counterfactual end-of-season ranking, but our rankings are in high agreement with respect to the teams that make the playoffs (95.65\%), the teams that get home court advantage in the playoffs (92.28\%), and the teams that receive double-digit lottery odds in the rookie draft (91.36\%).  We also provide a model extension that ensures each team's strength-of-schedule is not materially impacted by our choice of shortened season.  Finally, we conclude our paper in \cref{sec:conclusion}.




\section{Literature Review}\label{sec:lit-review}

A distinct feature 
of our study is the two-phase analytics approach that combines predictive and prescriptive models. We review existing literature in both streams: predicting single-game outcomes and end-of-season rankings, and scheduling sports leagues using algorithms and optimization.

\textbf{Predictive models for single-game outcomes.}
There is an extensive literature on predicting the outcome of a single sports game using historical data. On the one hand, there is an inherent difficulty in making such predictions, as the outcome of a  game depends both on luck and skill \citep{aoki2017luck}, and there is a limit to how much one can disentangle team/player skills from randomness \citep{martin2016exploring}. On the other hand, with ever-growing access to sports data and advancements in the fields of data analysis and machine learning, there has been a growing interest in predicting the outcomes of sporting events, both among researchers and for-profit organizations (e.g., \citealp{fivethirtyeight}). Existing models analyze the outcomes of sporting events using (i) Bayesian inference and rule-based reasoning \citep{miljkovic2010use},
(ii) Markov chain modeling \citep{kvam2006logistic,brown2010improved},
(iii) machine learning \citep{magel2014examining,prasetio2016predicting},
or (iv)
wisdom of crowds \citep{halberstadt1999effects}.
A key differentiator of our approach 
stems from the fact that we do not only predict the outcome of the next game to be played using all historical data available prior to that game. Instead, we predict the outcomes of all post-suspension games using only data from the pre-suspension period.

\textbf{Predictive models for end-of-season rankings.}
A few researchers have developed models to directly predict end-of-season team rankings given historical data of game play up to a certain (e.g., suspension) date. To the best of our knowledge, the first such paper is \cite{van2015predicting}, who study the final league rankings in European football leagues both before the start of the season and during the course of the season.
A body of research also focuses on determining the true ranking of individuals or teams in competitive sports based on network-based ranking systems \citep{motegi2012network,bozoki2016application}.
More recently, inspired by the COVID-19 suspension in the European football leagues, researchers \cite{van2021probabilistic} and \cite{csato2021coronavirus} studied the suspended season problem (also known as ``incomplete round robin league'') with the aim of predicting the final team ranking without additional games being played (i.e., if put into practice, this model would declare the league champion directly without resuming the season or playing any additional games).
These studies use descriptive and predictive models based on historical data from all games played prior to the suspension, to obtain a measure of strength for each team in the league and a measure of toughness of schedule. Using these two measures, the authors produce a projected final ranking and evaluate the accuracy of their proposed ranking using Monte Carlo simulation.
A key differentiator of our approach 
is our substantial prescriptive component. These papers do not predict outcomes of individual games and instead directly predict end-of-season rankings, which could be used to declare league standings without playing any additional games. In contrast, our predictive model produces predictions for each game that we then feed into our prescriptive model to determine the subset of games to play in our prescribed shortened season.

\textbf{Prescriptive models for scheduling sports games.} 
Within the sports scheduling literature, there are articles that focus on a primary objective (e.g., broadcast TV or travel logistics, fairness considerations), and others that propose multi-objective models.
Among those that use optimization to schedule basketball games, \cite{bean1980reducing} consider travel costs and player fatigue as the main goals, \cite{weiss1986bias} studies schedule bias between the regular season and post-season, while \cite{westphal2014scheduling} focuses on venue availability and broadcasting considerations. To propose a schedule
for NCAA basketball games, \cite{nemhauser1998scheduling} and \cite{henz2001scheduling} apply integer programming and constraint programming, respectively. Other papers develop tailored algorithms, often based on graph theory, for scheduling basketball games, see \cite{briskorn2009branch}.
We suggest the survey paper \cite{rasmussen2008round} for an overview of round-robin scheduling.
For a comprehensive list of articles in the broader scope of analytical methodologies applied to sports, including optimization and probabilistic modeling, see \cite{fry2012introduction, fry2012introduction1}.
Mixed-integer programs have been used to schedule games in different leagues: \cite{fleurent1993allocating} in hockey, \cite{goossens2009scheduling} in soccer, \cite{jiaqi2019designing} in baseball, and \cite{cocchi2018scheduling} in volleyball.

In contrast to the existing literature that schedules an entire season from a blank slate, the problem we consider compresses the remainder of an already-started season by selecting to play a subset of previously-scheduled games. To the best of our knowledge, this problem has not been previously studied. A significant novelty in our approach is our objective function, which attempts to achieve rankings similar to those that would have resulted had the season been played in full; this motivates us to introduce several model components which are novel in the context of sports scheduling, including ranking-based objectives, related stochastic optimization models, and finally predictive models for estimating the parameters of our stochastic optimization models.




\section{Problem}\label{sec:problem_background}

In this section, we first describe an NBA season and explain why a team's end-of-season ranking is important. We then frame the conclusion of the season as a problem of selecting a subset of the remaining games, which leads us to introduce several ranking similarity metrics.

\subsection{Background}\label{subsec:background}

\begin{figure}[ht]
    \centering
    \includegraphics[width=0.85\textwidth]{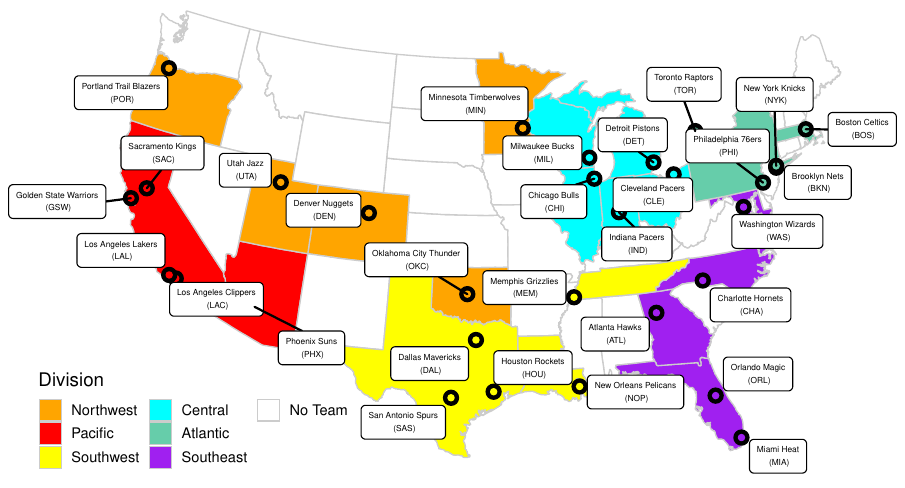}
    \caption{The Eastern Conference is comprised of the Central, Atlantic, and Southeast divisions, while the Western Conference consists of the Northwest, Pacific, and Southwest divisions.}
    \label{fig:nba-teams}
\end{figure}

The NBA is composed of 30 teams, divided into 2 conferences of 3 divisions with 5 teams each; for details, see Figure \ref{fig:nba-teams}.
In a regular season spanning approximately 180 days from October through April, each team plays 82 games according to the following formula: four games against the other four division opponents ($4 \times 4=16$ games), four games against six (out-of-division) conference opponents ($4 \times 6=24$ games), three games against the remaining four conference teams ($3 \times 4=12$ games), and finally two games against teams in the opposing conference ($2 \times 15=30$ games). A five-year rotation determines which out-of-division conference teams are played only three times. After 5 seasons, each team will have played 20 games against each in-division opponent, 18 games against each out-of-division opponent, and 10 games against each team from the opposing conference.

In our paper, we assume the regular-season standings directly determine the teams that advance to the playoffs. This follows the practice prior to the 2020--21 season, in which the 8 top-ranked teams in each conference (16 in the league) advance to the playoffs. However, starting with the 2020--21 season, only the top 6 teams in each conference automatically advance to the playoffs, with teams ranked 7 through 10 competing in a \emph{play-in tournament} to determine who is seeded $7^{\text{th}}$ and $8^{\text{th}}$ and who is eliminated; for details, see \cite{nba-play-in}. Note that the new format doesn't materially change our models or testing methodology, and since our problem instances are all from pre-COVID seasons, we exclude the play-in tournament from our analysis.

Once the standings are finalized, the playoffs begin. In each conference, the $i^{th}$-ranked team is initially matched to the $(9-i)^{th}$-ranked team, for $i \in \{1 \dots 8\}$. All playoff matchups are best-of-seven series, i.e., a team needs to win four out of seven games against the same opponent to win the matchup and progress to the next round (the loser of the matchup is eliminated). All matchups occur within-conference until the final matchup, which pits the winning team of the Eastern conference against the winning team of the Western conference. It is also important to note that the higher-ranked team in each matchup is awarded home court advantage; this means it hosts games 1, 2, 5, and 7, while the lower-ranked team hosts games 3, 4, and 6 (with games 5--7 played only if needed). Note that due to how teams are matched, the 4 top-ranked teams in each conference are given home court advantage in the first round of the playoffs.

There is a strong connection between a team's regular-season ranking and its playoff performance. Examining the history of 76 completed NBA seasons from 1946--2021, we find that (i) over 75\% of playoff series were won by the team with home court advantage, (ii) in only 5 seasons did an $8^{\text{th}}$-ranked team win a playoff series against a $1^{\text{st}}$-ranked team, and strikingly (iii) of the 76 NBA champions, 73 were ranked among the top 3 teams in the league at the end of the regular season.

To the extent that end-of-season rankings give teams preferential treatment in the playoffs which boost a team's chances of winning a championship, it is in the league's best interest to ensure that the ranking is fair, i.e., reflects to the greatest extent possible which teams are truly the best. Fairness can be in question when the season ends early. This is because when only a shortened season is played, it is possible for some teams to be matched with relatively easy-to-beat teams while others are matched with harder-to-beat teams, and this may result in a ranking that is quite different than one which would have resulted had the season been played in full. (We assume the full season's ranking is fair, since the league constructs the full season schedule in a balanced and equitable manner, and in general the public accepts the ranking at the end of the full season).

Finally, the order in which teams pick rookie players in the annual NBA draft is also tied to the final ranking, with the lowest-ranked teams given a higher chance of  drafting the best rookie player. Therefore, the ranking of teams outside of the top 8 in each conference is also important. The quality of the players in the draft varies from season to season, but some first-pick rookie players have included generational talents such as \emph{LeBron James}, \emph{Magic Johnson}, and \emph{Hakeem Olajuwon}, who have had huge impacts in leading their respective teams to win multiple championships.

\subsection{Problem Description}

\begin{figure}[t!]
    \centering
    \includegraphics[width=0.9\textwidth]{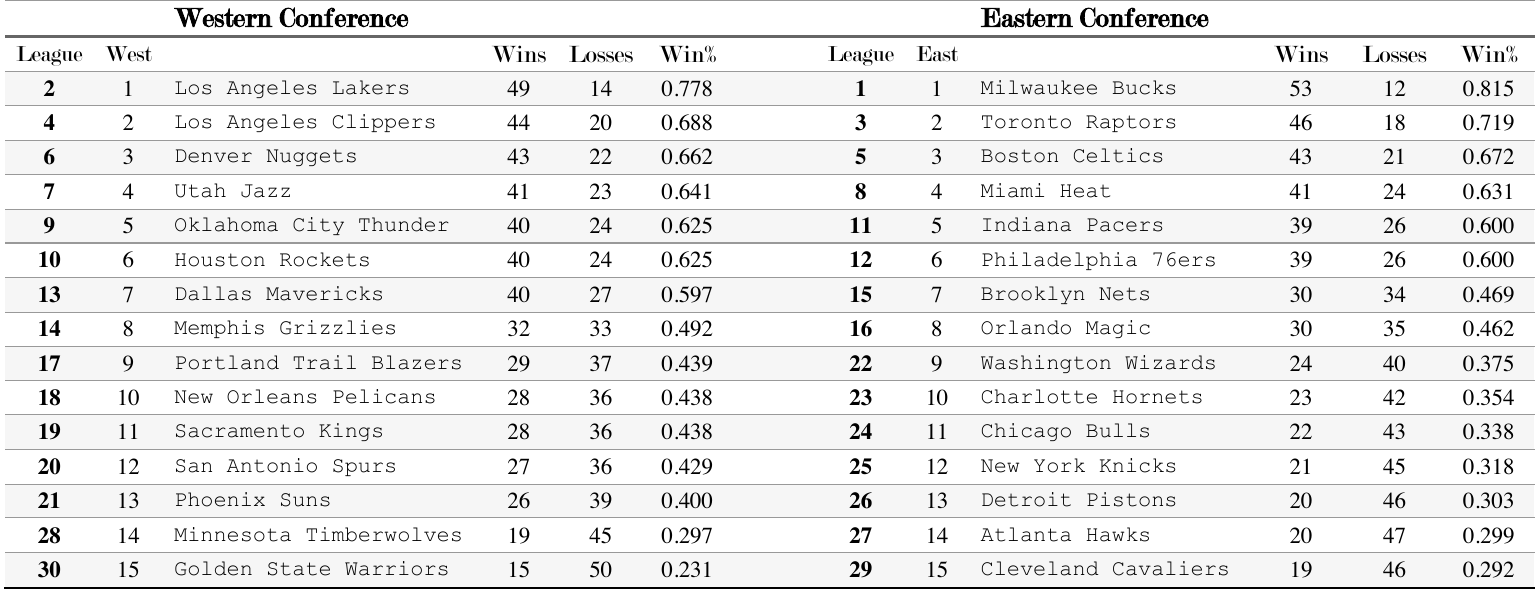} 
    \caption{NBA ranking at the time of suspension on March 11, 2020.\vspace{-0.5cm}}
    \label{fig:nba-ranking}
\end{figure}

At the time of the 2019--20 \mbox{COVID-19} suspension, 971 games in the 1230-game season were played leaving 259 games remaining; see Figure \ref{fig:nba-ranking} for the ranking at the time of suspension.
Given a target number of games each team should play in the full season, we wish to select a subset of the remaining 259 games that satisfies these targets. Typically, each of the 30 teams plays 41 home and 41 away games for a total of 82 games in the season. Shortening the season involves reducing the target from 82 games/team to a lower number (e.g., 70), with half the games at home and half away. Since the results of the 259 remaining games are uncertain, both the ranking produced by playing the full 82 game/team season and the ranking produced by playing a shortened (e.g., 70 game/team) season are uncertain. Our problem is to select a subset of games that minimizes the expected dissimilarity between the ranking of the full season and the ranking of the shortened season.  Before we introduce our models, we introduce several metrics that may be used for measuring the similarity of rankings.

\subsection{Measures of Similarity/Dissimilarity between Rankings}\label{sec:sim-dissim-measures}

We represent a ranking of $n$ teams as a vector, with components $1$ through $n$ permuted in order from the highest to the lowest percentage of games won during the regular season.  Throughout the paper, we follow the convention that $\hat{r}$ represents a ranking resulting from playing all games in the full season, while $r$ represents a ranking resulting from playing a specific subset of the remaining games (i.e., the shortened season case).  Furthermore, when we wish to distinguish between multiple rankings in the shortened season case, we use a superscript.  For example, $r^{(1)}$ and $r^{(2)}$ represent two distinct rankings resulting from concluding a shortened season with two different sets of games.

Two widely used measures of similarity/dissimilarity between rankings are \emph{Concordance}, and \emph{Euclidean distance}. We now define these metrics, using the following small example.

\textbf{Example 1:} Assume there are only four teams in the league: \texttt{LAL}, \texttt{MIL}, \texttt{LAC}, and \texttt{BOS}. Table \ref{table:ex1-rankings} contains the full-season ranking $\hat{r}$, and two alternative rankings $r^{(1)}$ and $r^{(2)}$.

\begin{table}[ht!]
\centering
\footnotesize
\begin{tabular}{cccc}
\toprule
Teams & Ranking ($\hat{r}$) & Ranking ($r^{(1)}$) & Ranking ({$r^{(2)}$}) \\ \midrule
\texttt{LAL}   & 1       & 1 & 4          \\ 
\texttt{BOS}   & 2       & 4 & 1            \\ 
\texttt{MIL}   & 3       & 2 & 3            \\ 
\texttt{LAC}   & 4       & 3 & 2           \\ \bottomrule
\end{tabular}
\vspace{5pt}
\caption{Example 1: Three different rankings.}
\label{table:ex1-rankings}
\end{table}
\vspace{-0.5cm}
\subsubsection{Concordance.}\label{sec:kendall-definition}
Concordance is a metric used to measure the ordinal association between two measured quantities, each with $n$ elements. Intuitively, concordance is high
when observations in two variables have similar ranks, and it is low
when observations have dissimilar (opposite) ranks.
Concordance, as a metric, is inspired by Kendall's rank correlation coefficient, or simply Kendall's $\tau$, introduced in \cite{kendall1938new}.
For a given pair of team rankings $(r, \hat{r})$, we call a pair of teams $(i,j)$ \emph{concordant} if either $i$ is above $j$ in both rankings (i.e., $\hat{r}_i > \hat{r}_j$ and $r_i > r_j$) or $i$ is below $j$ in both rankings (i.e., $\hat{r}_i < \hat{r}_j$ and $r_i < r_j$).
Conversely, we say a pair of teams $(i,j)$ is \emph{discordant} if their relative positions in the two rankings do not agree (i.e., either $\hat{r}_i > \hat{r}_j$ and $r_i < r_j$, or alternatively $\hat{r}_i < \hat{r}_j$ and $r_i > r_j$).
If $\hat{r}_i = \hat{r}_j$ or $r_i = r_j$, the pair of teams is neither concordant nor discordant.
Following our example, when we compare rankings $\hat{r}$ and $r^{(1)}$ from Table \ref{table:ex1-rankings}, the pair (\texttt{LAL}, \texttt{BOS}) is concordant, since in both rankings \texttt{LAL} stands higher than \texttt{BOS}. The pair (\texttt{BOS}, \texttt{MIL}) however is discordant, since \texttt{BOS} has the higher rank in $\hat{r}$, while \texttt{MIL} is higher in $r^{(1)}$.

Using the number of concordant pairs given two rankings, concordance ($\tau_C$) is defined as
\begin{equation}
    \tau_{C} = \text{number of concordant pairs}. \label{eq:kendall-tau-con-only}
\end{equation}
Note that $\tau_C$ is a number between 0 and $\binom{n}{2}$. Continuing our example, concordance is $\tau_C(r^{(1)},\hat{r}) = 4$ between $(r^{(1)}, \hat{r})$ and is $\tau_C(r^{(2)},\hat{r}) = 2$ between $(r^{(2)}, \hat{r})$.

\subsubsection{Euclidean distance.}\label{sec:spearman-rho-definition}
The (squared) Euclidean distance between two rank vectors is another metric used to measure dissimilarity between two alternative rankings, defined as
\begin{equation}
    d_{\text{E}}(r,\hat{r}) = \|r-\hat{r}\|^2 = \sum_{i \in T} (r_i - \hat{r}_i)^2,
    \label{eq:euclidean-dist}
\end{equation}
\noindent
where $T$ denotes the set of all the teams in the league. It can be easily verified that $d_{\text{E}}(\hat{r},r^{(1)})=6$ and $d_{\text{E}}(\hat{r},r^{(2)})=14$, which is consistent with our conclusion based on concordance $\tau_C$ that $r^{(1)}$, with shorter distance to $\hat{r}$, is the most similar to $\hat{r}$. We remark that Euclidean distance is equivalent to \emph{Spearman's Rank Correlation Coefficient} \citep{spearman1987proof} through an affine transformation with negative coefficient.
We also note the following relationship between concordance and Euclidean distance of two rankings. Proof of this proposition and other proofs are provided in Appendix~\ref{app:proofs}.
\begin{proposition}\label{prop:kendall_euclidean}
For arbitrary rankings $r$ and $\hat{r}$, the following relationship holds:
\begin{align}
    \sqrt{2 \,d_{\text{E}}(r,\hat{r})}\le n(n-1) - 2 \,\tau_{\text{C}}(r,\hat{r})
    \label{eq:kendall_euclidean}
\end{align}
\end{proposition}
\endproof

\subsubsection{Measuring Ranking Similarity/Dissimilarity Across Time.}
\begin{figure}[ht]
    \centering
	\subfigure[Concordance ($\tau_{\text{C}}(r,\hat{r})$)]{
		\includegraphics[clip,width=0.48\textwidth]{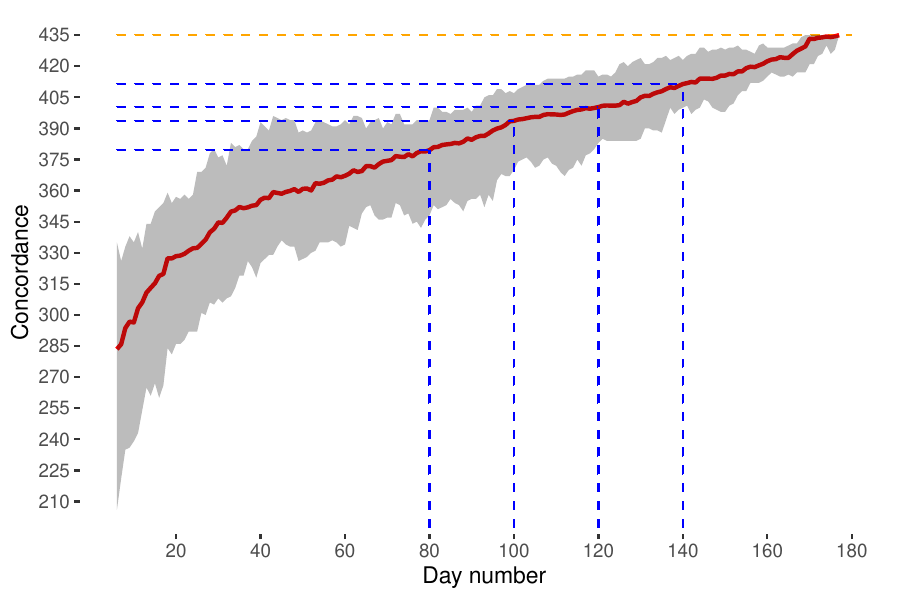}\label{fig:daily-concordance}}
	\subfigure[Euclidean distance ($\sqrt{d_{\text{E}}(r,\hat{r})}$)]{
		\includegraphics[clip,width=0.48\textwidth]{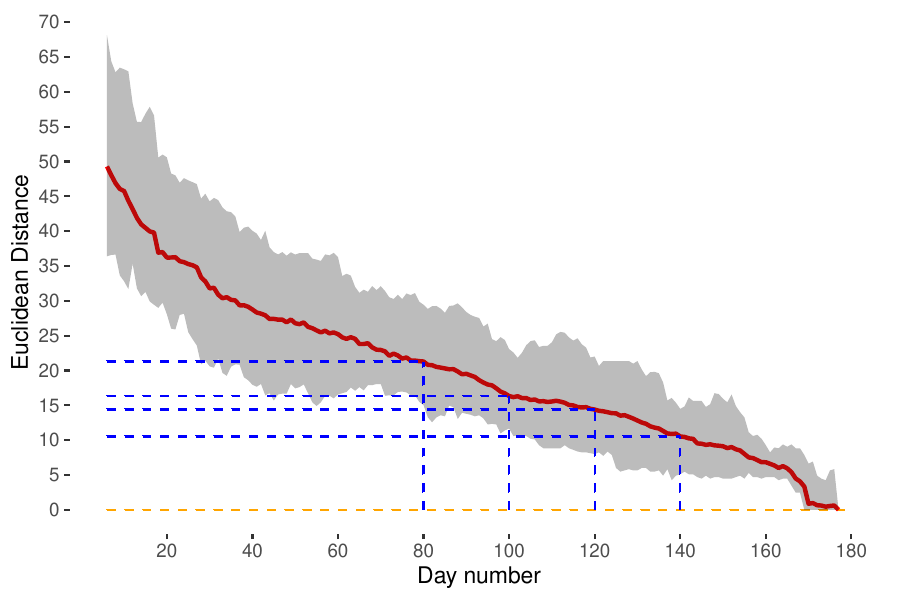}\label{fig:daily-euclidean}}
    \caption{
    Measuring ranking similarity/dissimilarity across time; grey area shows the range across 14 NBA seasons.}
    \label{fig:daily-metrics}
\end{figure}

It is interesting to plot how our metrics change over the course of an NBA season, as we compare daily team rankings ($r$) to the end-of-season ranking ($\hat{r}$). Note that with 30 teams in the NBA, the maximum concordance is $\binom{30}{2} = 435$. Figure \ref{fig:daily-concordance} plots concordance over time, while Figure \ref{fig:daily-euclidean} plots Euclidean distance. The red line is the mean over 14 NBA seasons from 2004--2018, while the grey band illustrates the range of values over these 14 seasons. As we approach the end of the season, concordance approaches its maximum of 435 while Euclidean distance converges to its minimum of 0.
\vspace{-0.5cm}




\section{Models}\label{sec:models}
\begin{figure}[ht]
    \centering
    \includegraphics[width=0.8\textwidth]{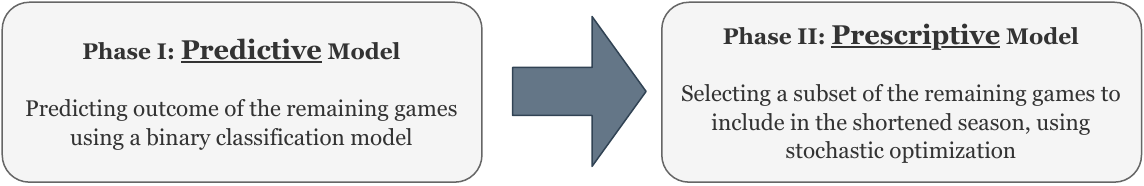}
    \caption{The two main phases of our methodology.}
    \label{fig:two-stages}
\end{figure}

Our two-phase modeling approach (c.f., Figure \ref{fig:two-stages}) falls under the ``predict-then-optimize'' paradigm \citep{mivsic2020data}. To choose the best subset of games for our shortened season, we need estimates of the outcomes of all remaining games.
Hence, in our first (``predictive") phase, we use historical data from all regular-season games played before the suspension to predict which teams win the remaining games. 
More specifically, we train binary classification models that incorporate game-related features (e.g., win percentage, point differential, and home-away indicator, among others). As we do not include team identifiers as categorical variables, we are able to take into account the toughness of the schedule throughout the season as experienced by each team.
\vspace{-0.5cm}
In our second (``prescriptive") phase, we determine which games to include and which games to cancel in the shortened season. 
We minimize the expected dissimilarity between the shortened season's ranking and the full season's ranking, where this expectation is taken over multiple possible scenarios that reflect the random chance each team has for winning each game. Specifically,
we treat game outcomes as Bernoulli random variables whose parameters are estimated in the first (``predictive") phase, and formulate our prescriptive models as stochastic optimization problems.

\subsection{Predictive Models}
\label{sec:predictive}

We now present our models for predicting the outcomes of games postponed due to the suspension. Since the response variable (i.e., whether the home team wins or loses) is binary, we model our prediction problem as a \textit{binary classification task}. 
More specifically, given a set of games $\scriptG$ from the pre-suspension period, each data point $g \in \scriptG$ used for training our model consists of a vector of $D$ features $\vec{x}_g=(x_{g,1},\dots,x_{g,D})$ and its binary class label $y_g$, where $y_g=1$ if the home team wins and $y_g=0$ otherwise. 
Our goal is to learn a discrimination rule $p:\mathbb{R}^D\to [0,1]$ representing the probability that a data point with feature vector $\vec{x}$ belongs to class 1 (``home team wins'').

We numerically evaluate 8 popular binary classification models from the literature. These include
Support Vector Machine \citep[SVM;][]{cortes1995support}, random forest \citep[RF;][]{breiman2001random}, bootstrap aggregating \citep[Bagging;][]{breiman1996bagging}, eXtreme Gradient Boosting \citep[XGBoost;][]{chen2016xgboost}, Extreme Learning Machines \citep[ELM;][]{huang2006extreme}, logistic regression (Logit), Gaussian Na\"ive Bayes (NB), and Multilayer Perceptron (MLP, also known as artificial neural network); see \cite{hastie2009elements} for descriptions of the last three. The explanatory variables in our dataset and specific features we use in our models are documented in \cref{sec:dataset-description} and Appendix \ref{app:exp-features}.

Most binary classifiers first estimate the class probability function $p(\vec{x})$, and then given each data point $\vec{x}$ use a threshold (e.g., $0.5$) to assign labels to observations (i.e., assign label 1 if $p(\vec{x}) > 0.5$, and 0 otherwise). One notable exception is traditional SVM, which directly assigns each data point to one of the two classes; here, to get probability estimates we follow common practice and employ a sigmoid calibration function~\citep{platt1999probabilistic,niculescu2005predicting}.

We perform model selection by evaluating our binary classification models using a performance metric, and choosing the best one.
Two applicable categories of performance metric are (1) those that compare predicted class (which is binary) with actual outcomes (also binary), and (2) those that compare predicted probability with actual outcomes.
The most natural criterion
in the first category is \emph{misclassification loss} (the complement of \emph{accuracy}, or the rate of correctly classified data points). 
But, since we are primarily interested in class membership probabilities rather than zero-one labels, 
we adopt a performance metric from the second category. Specifically, we use \emph{proper scoring rules}, which measure the quality of predicted probabilities~\citep{gneiting2007strictly}. 

For observations drawn from the distribution $F$, a scoring rule is called \textit{proper} if its expectation is maximized when the forecaster issues the probabilistic forecast $F$. Moreover, if $F$ is the unique maximizer, 
the scoring rule is \textit{strictly proper}.  
Examples of strictly proper scoring rules include the logarithmic, quadratic, and spherical scoring rules. Brier score and LogLoss, which measure the distance between estimated and true outcomes, are the complements of the quadratic and logarithmic scoring rules, respectively.
For details, see \cite{bickel2007some} and \cite{johnstone2011tailored}.

According to \cite{bickel2007some}, the logarithmic scoring rule favors \emph{sharper} probability values (i.e., those that are closer to zero or one). As suggested by \cite{johnstone2011tailored}, sharper probability values are preferred when the end use of estimated probabilities is a stochastic optimization problem, as we have in our second phase. This is because 
they reduce variance by focusing on the most likely scenarios. Hence, we use the logarithmic scoring rule as our model selection criterion to select the best-performing of our binary classifiers. In line with machine learning terminology, we henceforth use LogLoss defined below, which is the negative of the logarithmic scoring rule:
\begin{equation}
    \text{LogLoss} = -\frac{1}{|\scriptG|} \sum\nolimits_{g \in \scriptG} \Big(y_g \,\, \log(p_g) + (1 - y_g) \,\, \log(1 - p_g) \Big), \label{eq:logloss-def}
\end{equation}
where $\{p_g\}_{g\in \scriptG}$ is the vector of predicted probabilities and $\{y_g\}_{g\in \scriptG}$ is the vector of true outcomes. We compare the predictive performance of candidate models according to LogLoss in \cref{sec:pred-results}.

\subsection{Prescriptive Models}\label{sec:prescriptive}

In a league with $n$ teams, let $T$ denote the set of teams. Assume that at the time of suspension, a set $G$ of regular-season games remain to be played, and each team $i\in T$ has won a total of $\vwp^0_i$ games before the suspension. We 
represent each game $g\in G$ with a tuple $g=(i,j,k)$, where $i(g)\in T$ and $j(g)\in T$ denote the host and guest teams, respectively, and $k(g)$ indexes the $k^{\text{th}}$ match between these two teams (recall two teams may play each other more than once). We also define $G^h_i\subset G$ and $ G^a_i\subset G$ as the set of remaining home and away games for team $i$, respectively.

We model the outcome of game $g$ using the Bernoulli random variable $W_g$, which is one if the host team $i(g)$ wins, and zero if the guest team $j(g)$ wins.  For each game $g$, we estimate the parameter $p_g = P(W_g=1)$ using historical data as discussed in \cref{sec:predictive}. Formally, we denote the set of all possible outcomes of all games in the remainder of the full season by $\Xi$, and use $\xi \in \Xi$ to index a specific realization of all games' outcomes.  When helpful, we explicitly write $W_g(\xi)$ to indicate $W_g$'s dependence on $\xi$. 
For a given outcome $\xi\in\Xi$, the win percentage of team $i$ (total wins divided by games played) after playing all remaining games (i.e., at the end of the full regular season) is
\begin{align}
    \hat{\vwp}_i(\xi)=\frac{1}{\hat{m}}\left(y^0_i+\sum\nolimits_{g\in G^h_i}W_g(\xi)+ \sum\nolimits_{g\in G^a_i}(1-W_g(\xi))\right),
\end{align}
where 
$\hat{m}$ is the number of games played by each team in the full season (e.g., 82 for the NBA).
We will continue to use the caret ($\, \hat{} \,$) to denote quantities that correspond to the full regular season.

For each game $g\in G$, we define a binary decision variable $x_g$ that we set to one if game $g$ is included in the shortened season, and zero otherwise.
Note that these $x$-decisions are made before knowing the realization of $\xi$. 
We define $X$ as the set of feasible solutions, i.e., restrictions placed on the $x$-variables to express tactical/fairness considerations such as each team having the same number of home/away games in the season, as well as binary requirements on the $x$-variables, i.e.,
\begin{align}
  X = \left\{
  \begin{array}{@{}ll@{}}
    \sum_{g\in G^h_i}x_g=m^h_i, & \quad\forall i\in T \\
    \sum_{g\in G^a_i}x_g=m^a_i, & \quad\forall i\in T \\
    x_g\in\{0,1\}, & \quad\forall g\in G
  \end{array}\right\},\label{eq:home-away-constraints}
\end{align} 
where $m^h_i$ and $m^a_i$ denote the target number of home and away games for team $i$, respectively. For instance, if team $i$ has played 33 home and 31 away games so far before the suspension, and we decide to conclude the season with a total of 72 games for each team, then this team must play an additional $m^h_i=\frac{72}{2}-33=3$ home and $m^a_i=\frac{72}{2}-31=5$ away games. Another alternative would be to combine the constraints on the number of home/away games for each team into a single constraint of the form $\sum_{g\in G^h_i\cup G^a_i}x_g=m^h_i+m^a_i$ that sets a target for the total number of games to play without specific home/away sub-targets.

For a given shortened season $\xx\in X$ and realization $\xi\in\Xi$, we denote the win percentage of team $i$ at the end of the shortened season as $\vwp_i(\xx,\xi)$, where
\begin{align}
  \vwp_i(\xx,\xi)=\frac{1}{m}\left(y^0_i + \sum\nolimits_{g\in G^h_i}W_g(\xi) x_g+ \sum\nolimits_{g\in G^a_i}(1-W_g(\xi)) x_g\right),
  \label{eq:yi-def}
\end{align}
and $m$ is the target total number of games for each team in  the shortened season (e.g., 72).

Let $d(\vwp(\xx,\xi),\hat{\vwp}(\xi))$ be a measure of dissimilarity between the vectors $\vwp(\xx,\xi)$ and $\hat{\vwp}(\xi)$, i.e., a measure that compares the win percentage of each team at the end of the full season with the win percentage of each team at the end of the shortened season, for a specific shortened season $\xx$ and outcome $\xi$.
Note that there is a one-to-one correspondence between $\hat{\vwp}(\xi)$ and the team rankings at the end of the full season, and between $\vwp(\xx,\xi)$ and team rankings at the end of the shortened season. Therefore, $d(\vwp(\xx,\xi),\hat{\vwp}(\xi))$ can also be viewed as a dissimilarity measure between these rankings, and our goal is to find a shortened season $\xx$ that minimizes the expected value of this dissimilarity. 
That is, we are interested in solving stochastic optimization problems of the form
\begin{align}
    \min_{\xx\in X}\quad &\mathbb{E}_{\xi}[d(\vwp(\xx,\xi), \hat{\vwp}(\xi))],\label{opt:abstract}
\end{align}
for different choices of the dissimilarity measure $d$.  We now introduce two such formulations.

\subsubsection{Maximizing concordance.}\label{sec:model-pc}
For a given outcome $\xi$, let $\hat{r}(\xi)$ and $r(\xx,\xi)$ denote the ranking vector we get when the full season is played, and respectively, when the shortened season $\xx$ is played. We solve the following stochastic optimization problem to maximize the expected similarity between these two rankings using the average concordance metric:
\begin{align}
    \max_{\xx\in X}\quad &\mathbb{E}_{\xi}\left[\tau_{C}\left(r(\xx,\xi), \hat{r}(\xi)\right)\right].\label{opt:kentall-tau-abstract}
\end{align}

While this formulation is compact, its objective function is highly nonlinear; consequently, we linearize it as follows.  First, we define a parameter $\hat{z}_{ij}(\xi)$ which takes value one if team $i$ is above team $j$ in the full-season ranking $\hat{r}(\xi)$, and zero otherwise. Similarly, we introduce a binary variable $z_{ij}(\xx,\xi)$ which takes value one if team $i$ is above team $j$ in the shortened-season ranking $r(\xx,\xi)$, and zero otherwise.  Since $z_{ij}(\xx,\xi)+z_{ji}(\xx,\xi)=1$, we introduce only the $z_{ij}$-variables where $i<j$ and use $1-z_{ij}(\xx,\xi)$ in place of $z_{ji}(\xx,\xi)$ whenever it is needed.  As well, we introduce continuous variables $\vwp_i(\xx,\xi)$, $i \in T$, to keep track of the win percentage of team $i$ in the shortened season $\xx$ under realization $\xi$. Finally, since it is clear that the solution to this optimization problem encodes a single $x$-vector, we henceforth suppress the $x$-argument for the $\vwp$- and $z$-variables, and restate problem \eqref{opt:kentall-tau-abstract} as the following stochastic Mixed Integer Linear Program (MILP):
\begin{align}
    \mbox{[PC]}\;\max\quad & \;\mathbb{E}_{\xi}\left[\sum_{i\in T}\sum_{j\in T: j>i}\left( z_{ij}(\xi)\hat{z}_{ij}(\xi)+(1- z_{ij}(\xi))(1-\hat{z}_{ij}(\xi))\right)\right]\label{opt:kentall-tau-obj} \\
    \mbox{s.t.}\quad&\vwp_i(\xi)= \frac{1}{m}\left(y^0_i + \sum_{g\in G^h_i}W_g(\xi) x_g+ \sum_{g\in G^a_i}(1-W_g(\xi)) x_g \right) & \forall i\in T,\forall \xi\in\Xi \label{opt:kentall-tau-yi-def}\\
    &z_{ij}(\xi)-1 \le \vwp_i(\xi) - \vwp_j(\xi) \le z_{ij}(\xi) & \forall i,j\in T: i<j,\forall \xi\in\Xi \label{opt:kentall-tau-zij-def}\\
    &z_{ij}(\xi)\in \{0,1\} & \forall i,j\in T: i<j,\forall \xi\in\Xi\label{opt:kentall-tau-zij-domain}\\
    & \xx\in X.\label{opt:kentall-tau-x-domain}
\end{align}
The objective function \eqref{opt:kentall-tau-obj} counts the expected number of concordant pairs.
Constraint \eqref{opt:kentall-tau-yi-def} computes the win percentage for each team under each realization as defined by equation \eqref{eq:yi-def}, and constraints \eqref{opt:kentall-tau-zij-def} establish the relationship between win percentages and relative positions of teams.
\subsubsection{Minimizing Euclidean distance between win percentages.}\label{sec:MinWinPctDist}
We now consider the Euclidean distance between win percentages at the end of the shortened season $\vwp(\xx,\xi)$ and win percentages at the end of the full season $\hat{\vwp}(\xi)$. To minimize this dissimilarity measure, we  solve the following stochastic mixed integer quadratic program:
\begin{align}
    \mbox{[PW]}\;\min\quad &\mathbb{E}_{\xi}\left[\sum\nolimits_{i\in T}\left(\vwp_i(\xi)-\hat{\vwp}_i(\xi)\right)^2\right]\label{opt:win-percentage-diff-obj} \\
    \mbox{s.t.}\quad&\vwp_i(\xi)=\frac{1}{m}\left(y^0_i + \sum_{g\in G^h_i}W_g(\xi) x_g+ \sum_{g\in G^a_i}(1-W_g(\xi)) x_g\right) & \forall i\in T,\forall \xi\in\Xi \label{opt:win-percentage-diff-yi-def}\\
    & \xx\in X.\label{opt:win-percentage-diff-X}
\end{align}

While PW does not directly measure Euclidean distance between rankings, which is more closely tied to league outcomes than win percentages (see \cref{sec:spearman-rho-definition}), it has several computational advantages.
First, PW does not require binary variables $z_{ij}(\xi)$ and associated linking constraints, making it a lighter formulation than PC.
Moreover, as we shall show in \cref{sec:equiv-det-formulation},
we may derive a closed-form expression for the expected value in the objective function \eqref{opt:win-percentage-diff-obj}, which results in a much simpler deterministic equivalent problem, despite the objective being quadratic rather than linear.

\begin{proposition}\label{prop:euclidean_rank_winpercentage} Let $L$ be the least common multiple of $m$ and $\hat{m}$. There exists $D\le \frac{n(n^2-1)}{3} L^2$ such that
$$d_{\text{E}}(r(\xx,\xi),\hat{r}(\xi)) \le D \sum_{i\in T}\left(\vwp_i(\xx,\xi)-\hat{\vwp}_i(\xi)\right)^2 \quad \forall \xx\in X, \forall \xi\in \Xi$$
\end{proposition}
\endproof

This proposition formally shows that PW effectively minimizes the expected Euclidean distance between rankings. Obviously, the opposite does not necessarily hold (i.e., we can have identical rankings but different win percentages).




\section{Solution Methodology}\label{sec:solution-methods}

The stochastic optimization problems introduced in \cref{sec:prescriptive} contain $2^{|G|}$ realizations of $\xi$, each with their own set of second-stage decision variables and constraints. As the full stochastic optimization problems are too large to solve directly, we introduce an exact deterministic reformulation paired with a fast tailored algorithm for approximately solving PW.
Next, we develop a robust optimization reformulation of PW. Finally, we also provide two solution methods for approximately solving PC which replace the random variables with (i) their means and (ii) a small finite number of samples (see Appendix \ref{app:approx-solution-methods-pm-PC} for these methods, which we will refer to as PC-MVP and PC-SAA, respectively, in the computational experiments).

\subsection{Equivalent Deterministic Reformulation}\label{sec:equiv-det-formulation}

We now show how the stochastic problem PW from \S\ref{sec:MinWinPctDist} can be solved using an equivalent deterministic problem. We will use the notation $\mathbb{V}$ to refer to the variance of a random variable.

\begin{theorem}\label{theorem-pw}
The stochastic model PW can be solved using the following equivalent deterministic linearly-constrained convex quadratic mixed-integer optimization problem:
\begin{align}
 \emph{\mbox{[PW-DQIP]}}\;\min\quad &\sum\nolimits_{i\in T}\left(\left(\mu_i-\hat{\mu}_i\right)^2 +v_i\left(1-\frac{2m}{\hat{m}}\right)+ \hat{v}_i \right) \label{opt:win-percentage-diff-det-obj}\\
   \text{s.t.}\quad &\mu_i = \frac{1}{m}\left(\vwp^0_i + \sum\nolimits_{g\in G^h_i}p_g x_g+ \sum\nolimits_{g\in G^a_i}(1-p_g) x_g\right) & \forall i \in T\label{opt:win-percentage-diff-det-mu}\\
    &v_i = \frac{1}{m^2}\sum\nolimits_{g\in G^h_i\cup G^a_i}p_g(1-p_g) x_g & \forall i \in T \label{opt:win-percentage-diff-det-v}\\
    & \xx\in X,\label{opt:win-percentage-diff-det-X}
\end{align}
where the decision variables, in addition to $\xx = \{x_g, g \in G\}$, include $\mu_i$ and $v_i$ which encode the mean and variance, respectively, of the win percentage of team $i$ in the shortened season. Moreover, the following parameters represent the mean and variance, respectively, of the win percentage of team $i$ in the full season: 
\begin{align}
    \hat{\mu}_i= \mathbb{E}_{\xi}\left[\hat{\vwp}_i(\xi)\right]= &\frac{1}{\hat{m}}\left(\vwp^0_i + \sum\nolimits_{g\in G^h_i}p_g + \sum\nolimits_{g\in G^a_i}(1-p_g)\right) \label{opt:win-percentage-diff-det-muHat} \\
    \hat{v}_i= \mathbb{V}_{\xi}\left[\hat{\vwp}_i(\xi)\right]=&\frac{1}{\hat{m}^2}\sum\nolimits_{g\in G^h_i\cup G^a_i}p_g(1-p_g). \label{opt:win-percentage-diff-det-vHat}
\end{align}
\end{theorem}

\subsection{Frank-Wolfe Algorithm}\label{sec:frank-wolfe-algorithm}
For the purpose of developing an efficient way to solve PW-DQIP, we now turn our attention to the combinatorial structure of PW-DQIP's feasible region. 

\begin{proposition}\label{prop:tum}
The coefficient matrix of the set of feasible schedules $X$ is totally unimodular.
\end{proposition}

As a result of Proposition~\ref{prop:tum}, and given that the right-hand-side values in $X$ (i.e., $\{m_i^h\}$ and $\{m_j^a\}$) are integral, optimizing a linear function over the continuous relaxation of $X$ (denoted $\bar{X}$) using the Simplex method yields an integral optimal solution.
This property of PW-DQIP lends itself well to the Frank-Wolfe (FW) method \citep{frank1956algorithm, jaggi2013revisiting}. FW is an algorithm for solving non-linear convex optimization problems of the form 
$\min_{\xx\in \bar{X}}f(\xx),$
where $f$ is a smooth convex function and $\bar{X}$ is a compact convex set. 
At each iteration $t$,
FW replaces $f$ with its linear approximation at an incumbent point $\xx^{(t)}\in \bar{X}$, to produce an ``atomic'' solution
\begin{align}
    \hat{\xx}^{(t)}=\argmin_{\xx\in \bar{X}} \nabla f(\xx^{(t)})^{\top}\xx, \label{fw-sub}
\end{align}
and then performs a line search between $\xx^{(t)}$ and $\hat{\xx}^{(t)}$ to produce the next iterate $\xx^{(t+1)}\in \bar{X}$.

Algorithm~\ref{pseudo-code-frank-wolfe} in Appendix \ref{app:frank-wolfe-algorithm} presents our implementation of the FW algorithm for solving the continuous relaxation of PW-DQIP. 
Our FW implementation is particularly efficient, since $f$ is a convex quadratic function, and so (i) its gradient is easily computed, and (ii) a closed-form optimal solution to the line search step can be found using the first-order optimality conditions.
Moreover, given that the atomic solution $\hat{\xx}^{(t)}$ produced by solving the transportation problem \eqref{fw-sub} is a feasible integer solution, it provides an upper bound on the optimal value of PW-DQIP. As the FW algorithm iterates, ${\xx}^{(t)}$ converges to the optimal fractional solution and the upper bound $\hat{\xx}^{(t)}$ becomes progressively tighter. Upon terminating at a finite iteration $t$, we use the best integer-feasible point $\hat{\xx}^{(t)}$ found thus far, over iterations $1\dots t$, as a near-optimal integer solution to PW-DQIP. Henceforth, we will use PW-FW to refer to producing a PW solution using FW.

\subsection{Robust Optimization Reformulation}\label{sec:robust-opt-formulation}

In the data science literature, there are many examples of successful ensemble models, which combine the results of multiple predictive models to produce more robust results \citep[c.f.,][]{sagi2018ensemble}.
With this in mind, we derive a Min-Max Regret (MMR) formulation that combines the predicted probabilities from several of our predictive models, allowing our prescriptive model's results to be more robust to misspecification and model overfitting.

Given a set $L$ of candidate predictions indexed by $l$, we define the following parameters: (i) $p^{(l)}_g$ is the probability that the home team wins game $g$ under candidate $l$; (ii) $\theta^{(l)}$ is the optimal value (or, alternatively, a lower bound on the optimal value) of PW-DQIP when solved using candidate $l$; and (iii) $\hat{\mu}_i^{(l)}$ and $\hat{v}_i^{(l)}$ correspond to the parameters defined by \eqref{opt:win-percentage-diff-det-muHat} and \eqref{opt:win-percentage-diff-det-vHat}, respectively, when $p^{(l)}_g$ is used in place of $p_g$. In addition to $x_g$, the MMR counterpart of our PW model has decision variables $\mu_i^{(l)}$ and $v_i^{(l)}$, which are the candidate-specific versions of $\mu_i$ and $v_i$ from \eqref{opt:win-percentage-diff-det-mu} and \eqref{opt:win-percentage-diff-det-v}, respectively. Finally, $\theta$ is a decision variable that captures the maximum regret over using all candidate win-probability vectors in the PW objective. Minimizing this maximum regret yields the following convex mixed-integer quadratically-constrained program:
\begin{align}
\text{[PW-MMR]}\;\min\quad & \theta \label{opt:win-percentage-diff-MMR-obj1}\\
\text{s.t.}\quad &\theta\ge \sum\nolimits_{i\in T}\left(\left(\mu^{(l)}_i-\hat{\mu}^{(l)}_i\right)^2 +v^{(l)}_i\left(1-\frac{2m}{\hat{m}}\right)+ \hat{v}^{(l)}_i \right) - \theta^{(l)} & \forall l \in L \label{opt:win-percentage-diff-MMR-obj2}\\
   &\mu^{(l)}_i = \frac{1}{m}\left(\vwp^0_i + \sum\nolimits_{g\in G^h_i}p^{(l)}_g x_g+ \sum\nolimits_{g\in G^a_i}(1-p^{(l)}_g) x_g\right) & \forall i \in T, \; l \in L \label{opt:win-percentage-diff-DRO-mu}\\
    &v^{(l)}_i = \frac{1}{m^2}\sum\nolimits_{g\in G^h_i\cup G^a_i}p^{(l)}_g(1-p^{(l)}_g) x_g & \forall i \in T, \; l \in L \label{opt:win-percentage-diff-DRO-v}\\
    & \xx\in X.\label{opt:win-percentage-diff-DRO-X}
\end{align}




\section{Computational Experiments}\label{sec:comp-exp}
In this section, we (1) evaluate our predictive models and choose the best one(s) to use in our prescriptive phase (``model selection"), and (2) solve all variants of our prescriptive models using the home-team win-probabilities estimated during our predictive phase, and evaluate the quality of the shortened seasons produced by our prescriptive models.
All models (predictive and prescriptive) were coded in \texttt{Python 3.8.8}. For our predictive models, we used the \texttt{scikit-learn} package \citep{scikit-learn}, and for solving our mixed integer programs in our prescriptive models we used \texttt{Gurobi 9.5.0} with all solver settings left at their default values.
All experiments were conducted on a computer with a 2.6 GHz Intel Core i7 CPU and 16 GB of memory.

\subsection{Dataset Description}\label{sec:dataset-description}
We use historical data from 14 NBA seasons (2004--2010, 2012--2018), which are the years with the same regular season structure as today; i.e.,  30 teams, each playing 82 games with schedules constructed in the manner described in \cref{subsec:background}. We omit the shortened seasons 2010--11 and 2020--21, and do not consider data prior to 2004 as back then the NBA had fewer teams. 
We used the box score datasets publicly available on the NBA's official website \citep{nba-stats} which contains detailed information for each game, team and player. From this, we created 56 datasets that in turn consider, for each of the 14 seasons, 4 alternative hypothetical suspension dates (i.e., days 80, 100, 120, 140 of the season). Note that each regular season spans between 170--180 days.

In our experiments, we fixed the number of games per team in the shortened season so approximately half the post-suspension games are scheduled and the other half cancelled (alternative targets can easily be achieved by adjusting the parameters $m_i^h$ and $m_i^a$, but we felt that varying the suspension day in our experiments already provides sufficient sensitivity analysis).
Table \ref{table:suspension_instances} provides some summary statistics of our instances. For each suspension day, we list the number of Games per Team ($GT$) we wish to have in the shortened season, the average number of Games Played ($GP$) per team prior to suspension, the average number of Games we will Schedule ($GS$) per team post-suspension, and the number of Games we will Cancel ($GC$) per team post-suspension. Note that (i) $GP+GS+GC=82$ since a full season has 82 games per team; and (ii) $GP+GS=GT$.

\begin{table}[ht!]
\centering
\footnotesize
\renewcommand{\arraystretch}{1.2}
\scalebox{1}{
\begin{tabular}{@{}ccccc@{}}
\toprule
Suspension & GT & GP & GS & GC \\ \midrule
Day 80\;\; & 62     & 38.4 & 23.6 & 20 \\
Day 100    & 66     & 48.6 & 17.4 & 16 \\
Day 120    & 70     & 56.4 & 13.6 & 12 \\
Day 140    & 74     & 66.0 & 8.0  & 8  \\ \bottomrule
\end{tabular}}
\vspace{0.2cm}
\caption{Summary statistics of our problem instances, averaged over all teams and all 14 NBA seasons.}
\label{table:suspension_instances}
\end{table}

\subsection{Predictive Model Results} \label{sec:pred-results}
In this subsection, we describe the explanatory features used in our predictive models, the pre-processing we performed on our datasets, and the cross-validation we used to assess our predictive models. Finally, we compare the predictive performance of the 8 binary classifiers described in \cref{sec:predictive}.

\subsubsection{Experimental setup.}\label{sec:pred-model-exp-setup}
We use the basic and advanced statistics published in the box score datasets ~\citep{nba-stats} to curate four groups of explanatory features, including i) overall team performance, ii) basic team-level statistics, iii) advanced team-level statistics, and iv) player-level statistics.
When taken together, this results in 128 features, the details of which we 
provide in Appendix \ref{app:exp-features}.
In our independent exploratory analysis, we have found that the features that tend to be the most important are, ranked in order from most to least important: i) overall team performance measured primarily by win percentage, ii) style of play metrics such as pace and number of possessions created, iii) shooting efficiency, iv) offensive fire power measured by metrics such as effective field goal percentage, true shooting percentage, points in the paint, and offensive rating in general, v) rebounding, and finally vi) assists and ball movement.
We first normalize all features to the $[0,1]$ range, and then use \emph{Principal Component Analysis (PCA)} to eliminate multicollinearity and prevent overfitting.
We retain only the first 25 principal components with highest eigenvalue, which explain more than 90\% of the total variance in our training data.

We use 5-fold cross-validation to tune our predictive models' hyper-parameters.
For this, we partition our data into \emph{training, validation,} and \emph{test} datasets. Pre-suspension games are divided into \emph{training} and \emph{validation} sets through 5 folds at random, with the validation dataset having the same size ($30\%$ of the of pre-suspension dataset) across all folds. Post-suspension games comprise the \emph{test} dataset. On each fold, we fit each classifier using a training dataset and evaluate its performance using the corresponding validation dataset. Then, we measure the performance of each predictive model using the average Logloss across the five folds.
We reserve the test datasets for evaluating the combined performance of our predictive and prescriptive models.

As is well-known in machine learning, not all classification models produce unbiased class probabilities. Therefore, we employ Platt scaling \citep{platt1999probabilistic}, which involves transforming the home team's win probabilities using a sigmoid calibration function to recover unbiased win probability estimates.
For this, we use the Python package \texttt{scikit-learn}, which uses multiple randomly-generated held-out samples from the training dataset to compute the calibrated probabilities.

\subsubsection{Predictive performance.}\label{sec:pred-performance}
We evaluate the performance of our predictive models using LogLoss as defined in \eqref{eq:logloss-def}, which is a strictly proper scoring rule (see \cref{sec:predictive} for details). We also report the \textit{accuracy} and the \textit{Area Under the receiver operating characteristic Curve} (AUC), which are not proper scoring rules but are widely used for evaluating the performance of binary classifiers. Accuracy measures the proportion of correctly predicted class labels. AUC summarizes the trade-off between true positive rate and false positive rate when the threshold to predict class labels changes \citep[see e.g.,][]{fawcett2006introduction}. Larger AUC values generally correspond to better predictive performance.

The LogLoss metric plays a crucial role in three distinct elements of our predictive phase. First, it serves as the internal loss function for all but three of our binary classifiers, with notable exceptions being SVM, Na\"ive Bayes, and Random Forest, which utilize hinge loss, Gini impurity, and no specific loss function, respectively; see \citep{shen2005loss} for additional details. Second, LogLoss guides our hyper-parameter optimization process across all eight predictive models. Lastly, we use it as our model selection criterion to identify the most effective of our eight candidate models.

\begin{figure}[t!]
\centering
\includegraphics[width=\textwidth]{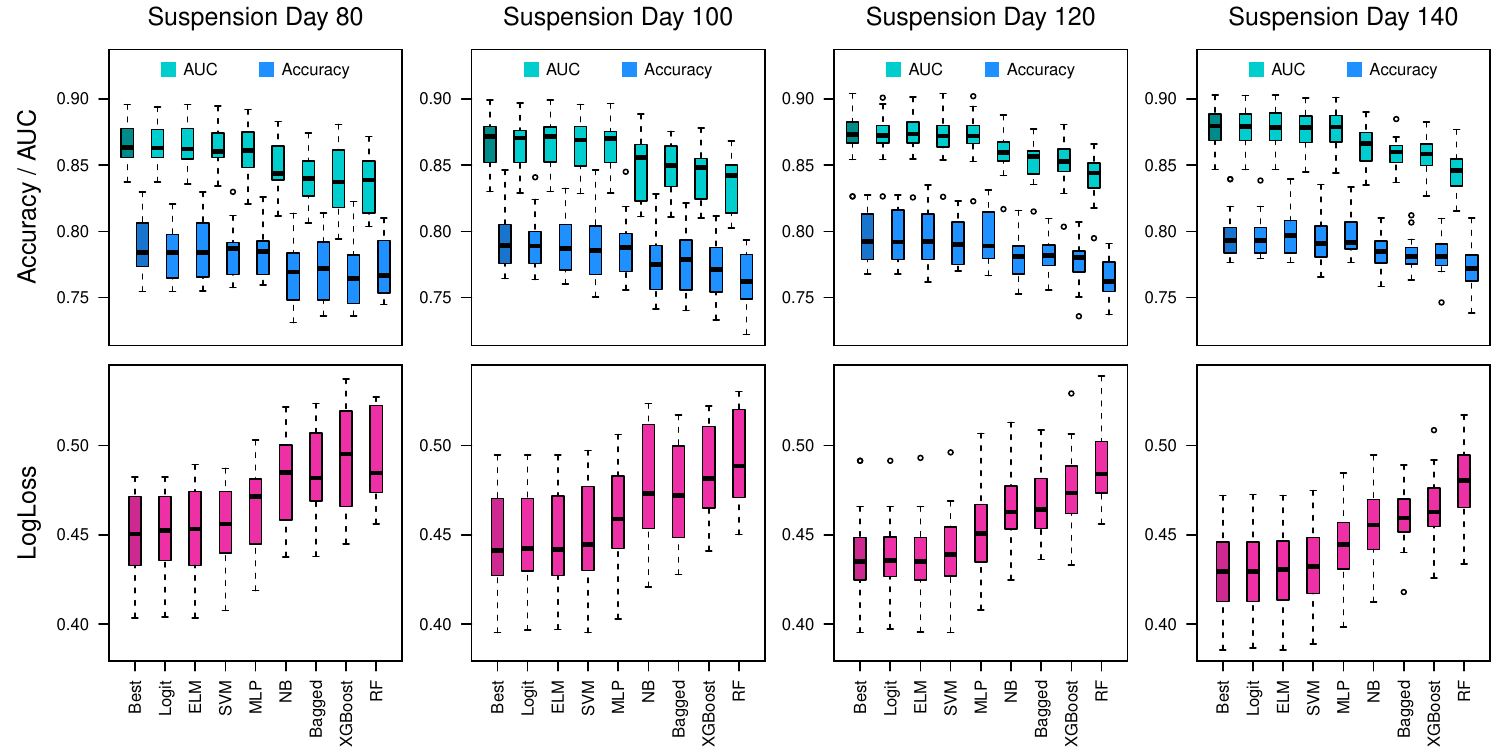}
\caption{Accuracy, AUC, and LogLoss for 9 different classifiers. The distributions are measured using the validation datasets across 14 NBA seasons.}
\label{fig:predictive_performance}
\end{figure}

Figure \ref{fig:predictive_performance} depicts, for 4 choices of suspension day and 9 classifiers, the distributions of accuracy, AUC, and LogLoss across 14 NBA seasons. In addition to the 8 classifiers described in \S\ref{sec:predictive}, we also define ``Best'' as a composite classifier that, for each season, employs the lowest-LogLoss classifier among the 8 classifiers tested. Each boxplot's underlying distribution corresponds to 14 values (one per season), where each season's performance is taken as an average of the model's performance over the five validation datasets (cross-validation folds). Box heights correspond to Inter-Quartile Ranges (IQR), with medians marked by dark horizontal lines inside each box.
Whiskers mark the most extreme-valued data points within 1.5IQR units above and below the boundary of the box, and any points lying outside this range are marked as outliers.

As seen in Figure \ref{fig:predictive_performance}, the Logit, ELM, SVM and MLP models generally have the lowest LogLoss values, which indicates their superiority in predicting well-calibrated and unbiased probabilities. The same four models also outperform the other models according to accuracy and AUC, with Logit being the better model overall. It is worth noting that our best predictive models achieve an accuracy over 75\%, which is in line with the accuracy of models that predict the outcome of the next game to be played (see, for example \cite{kvam2006logistic,brown2010improved}).
 
We conjecture that two factors drive the performance of our predictive models:
the level of league competition and how early we suspend the season. In a more competitive season, it becomes more difficult to predict game outcomes, since teams are more evenly-matched.
Second, as we increase the suspension day from 80 to 140, the pre-suspension dataset grows and as a result, our predictive models generally perform better in all three metrics depicted in Figure \ref{fig:predictive_performance}, both in terms of the average value as well as the variation across seasons (i.e., length of the boxplot).

In general, since different predictive models perform better in different seasons, we use our composite ``Best" classifier to produce the home-team win-probabilities for the post-suspension games used  by our prescriptive models. That is, for each season and suspension day, we select the lowest-LogLoss of the 8 classifiers as measured by 5-fold cross-validation on our validation dataset, and then train the selected classifiers on all pre-suspension games to produce home-team win-probabilities for all post-suspension games.


\subsection{Prescriptive Model Results}\label{sec:prescriptive-results}

We now evaluate our prescriptive models and provide practical insights that apply to their use. 

\subsubsection{Experimental setup.}
\label{sec:presc-model-exp-setup}

We measure the runtimes and solution quality of all prescriptive models discussed in \cref{sec:prescriptive} (i.e., PW-DQIP, PW-FW, and PW-MMR), as well as the two introduced in Appendix \ref{app:approx-solution-methods-pm-PC} (i.e., PC-MVP and PC-SAA); results are summarized in Table~\ref{table:prescriptive_models}.
As described in detail in  Appendix~\ref{app:saa_sample_size}, the trade-off between solution quality and runtime of SAA is balanced at 50 scenarios. Thus, we use 50 scenarios for our PC-SAA model. 
As well, to improve the computational efficiency of our PC model, we introduce variable fixing techniques (see Appendix \ref{app:variable-fixing}) to deduce and fix many $z$-variables at their optimal values. As detailed in Table~\ref{table:elimination} in Appendix \ref{app:variable-fixing}, our method eliminates 75\% of the $z$-variables in PC-MVP and 60\% in PC-SAA, which is significant, as it translates to closing the optimality gap at a faster rate (e.g., for suspension day 80, after 3,600 seconds of running PC-SAA the gap is $2.02\%$ with variable fixing, and $6.93\%$ without).

\begin{table}[t!]
    \centering
    \footnotesize
    \renewcommand{\arraystretch}{1.2}
    \scalebox{0.9}{
	\begin{tabular}{@{}lll@{}}
        \toprule
        \textbf{Model} & \textbf{Objective} & \textbf{Solution Method} \\ \midrule
        \textbf{PW-DQIP} & Min. distance of win percentages & Deterministic equivalent  Quadratic Program \\
        \textbf{PW-FW} & Min. distance of win percentages & Produce near-optimal solution for PW-DQIP using Frank-Wolfe \\
        \textbf{PW-MMR} & Min. distance of win percentages & A Robust reformulation of the PW-DQIP model \\ \midrule
        \textbf{PC-MVP} & Max. concordance of rankings & Replace random variable $W_g$ with its expected value $p_g$ for each game $g\in G$ \\
        \textbf{PC-SAA} & Max. concordance of rankings & Replace distribution $\Xi$ with sample $\mathcal{S}$, and expected value with sample average\\
         \bottomrule
    \end{tabular}}
    \vspace{0.2cm}
    \caption{Summary of prescriptive models.}
\label{table:prescriptive_models}
\end{table}

We implement two benchmarks. First, we are interested to know how well our models perform relative to an approach that does not explicitly optimize the end-of-season ranking when selecting the games in the shortened season. For this, we implement a greedy heuristic (henceforth known as ``Greedy") that selects games according to their original scheduled dates, with earlier games assigned first until the target number of games for each team are met. Second, we are also interested to know how much better we can do by playing our optimally-chosen shortened season relative to not playing any more games after the suspension (the ``Status Quo" ranking).

In line with previous studies in sports analytics (e.g., \citealp{van2021probabilistic, chater2021fixing,csato2021simulation,sziklai2022efficacy}), we evaluate the expected performance of our models using Monte Carlo simulation. More specifically, we draw the outcomes of each post-suspension game from a Bernoulli distribution with home-team win-probability estimated by our ``Best" predictive model.
To measure the expected performance,
we generate a sample of 10,000 game outcomes from these Bernoulli distributions. Each realization yields two rankings, one at the end of the shortened season, and the other at the end of the full season with all games played. We then compare these rankings using our primary performance metric, the number of concordant pairs between rankings.

Note that we use estimated home-team win-probabilities in both our prescriptive models and our simulation. To ensure robustness of our testing methodology, we use some \emph{unseen} data points to estimate the home-team win-probabilities for our simulation. To this end, we hold out 20\% of the pre-suspension data at random when we estimate the parameters ${p_g, \; g \in G}$, used as home-team win-probabilities in our prescriptive models, and use the entire pre-suspension dataset to estimate the home-team win-probabilities used by our simulation. 

Finally, as an additional robustness check, we also ran simulations with probabilities estimated using multiple evaluator models with different functional forms. For each season and suspension day, we compared PW-FW to Greedy using all combinations of prescriptive model parameters ${p_g, \; g \in G}$, and simulation Bernoulli probabilities estimated by our 8 classifiers from \S\ref{sec:predictive} (for a total of 64 combinations). As illustrated in Appendix~\ref{app:cross_sim}, all PW-FW solutions outperformed the Greedy solutions by a considerable margin.

\subsubsection{Prescriptive model results.}
\label{sec:presc-model-results}

Table \ref{table:prescriptive_performance_results} summarizes the results of running all four prescriptive models on the instances described in Table \ref{table:suspension_instances}, averaged over 14 seasons. The ``Time'' column reports runtimes in seconds, while ``Conc.'' is  concordance ($\tau_C$) measured by our simulation. PC-MVP and PW-FW solved all instances to optimality, while the optimality gap of PW-DQIP was never more than $10^{-3}$. For PC-SAA and PW-MMR, no instance converged to optimality by the 3600-second time limit. We report the percentage optimality gap for PC-SAA, while for PW-MMR we provide the absolute optimality gap ($UB - LB$). The latter metric is more appropriate for PW-MMR, since its lower bound is always zero. 

Comparing our prescriptive models, we observe the following. PC-MVP is fastest, but solution quality (concordance) is lowest. PC-SAA produces higher-quality solutions than PC-MVP, but is slow. All PW-based models generate high-quality solutions. While PW-DQIP typically produces slightly higher-quality solutions than PW-FW, our Frank-Wolfe algorithm is over \textbf{a thousand} times faster than the integer program of PW-DQIP. With regard to PW-MMR, it is interesting that although its solutions are more robust, overall solution quality is comparable to the other PW methods. While runtimes of PW-MMR are slower as we solve a mixed-integer quadratically-constrained program, it is worth pointing out that it may be possible to extend our Frank-Wolfe method to PW-MMR (see our comment at the end of Appendix \ref{app:frank-wolfe-algorithm}; we leave the details to future work).
Finally, we also observe that as the season is suspended later, there are fewer games to choose from, and all models perform better as the solution space shrinks.

\begin{table}[ht]
\centering
\footnotesize
\renewcommand{\arraystretch}{1.2}
\scalebox{1}{
\begin{tabular}{@{}rllllllllllll@{}}
\toprule
\textbf{Sus. Day}        & \multicolumn{2}{c}{PC-MVP} & \multicolumn{3}{c}{PC-SAA} & \multicolumn{2}{c}{PW-DQIP} & \multicolumn{2}{c}{PW-FW} & \multicolumn{3}{c}{PW-MMR} \\ \cline{2-13}
\multicolumn{1}{r}{(GT)} & Time        & Conc.        & Time   & Gap     & Conc.   & Time         & Conc.        & Time       & Conc.        & Time & Abs. Gap        & Conc.        \\ \midrule
80 (62)                  & 0.08        & 413.35       & 3600   & 4.98\%  & 414.24  & 3600         & 415.96       & 0.39       & 415.86       & 3600 & 0.00037       & 415.84       \\
100 (66)                 & 0.09        & 414.60       & 3600   & 1.66\%  & 415.80  & 3600         & 417.94       & 0.24       & 417.83       & 3600 & 0.0002       & 417.87       \\
120 (70)                 & 0.10        & 417.43       & 3600   & 1.49\%  & 418.55  & 934.11       & 419.96       & 0.12       & 419.88       & 3600 & 0.00015       & 419.92       \\
140 (74)                 & 0.06        & 420.64       & 3600   & 1.03\%  & 421.45  & 82.33        & 422.24       & 0.05       & 422.20       & 3600 & 0.00012        & 422.26 \\ \bottomrule     
\end{tabular}}
\vspace{0.2cm}
\caption{Performance of the prescriptive models, averaged over 14 seasons}
\label{table:prescriptive_performance_results}
\end{table}

\begin{figure}[t!]
    \centering
	\includegraphics[clip,width=\textwidth]{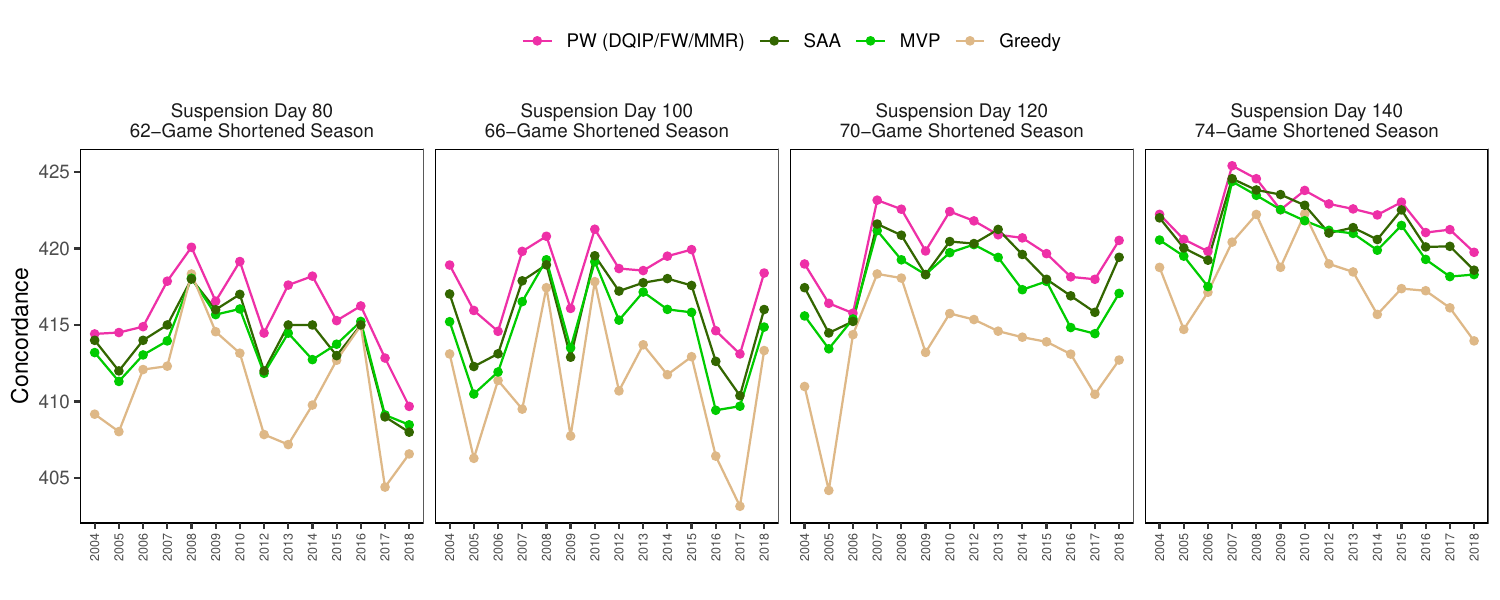}
    \caption{Simulation results for prescriptive models (concordance) across 14 NBA seasons.} 
    \label{fig:sim_results}
\end{figure}

As different seasons unfold in different ways, it is interesting to compare our prescriptive models on a per-season basis. Figure \ref{fig:sim_results} plots concordance as measured by our simulation for each of the 14 NBA seasons and 4 suspension days. The three PW models (i.e., PW-DQIP, PW-FW, and PW-MMR) have very similar solution quality and thus we combine their plots. From Figure \ref{fig:sim_results}, we see that all five proposed models outperform the benchmark greedy algorithm on all instances (i.e., for all seasons and suspension days). Typically, PW dominates PC-SAA and PC-MVP as well. 

We provide a few comments to explain the differences in performance across these models, which all approximate our stochastic optimization model PC in distinct ways. First, we note that although both PC-MVP and PC-SAA maximize concordance directly, PC-MVP uses only the means of the random variables while PC-SAA samples a small number of scenarios. Because PC-SAA more faithfully represents PC, it is not surprising that it typically outperforms PC-MVP. On the other hand, it is interesting to observe that the PW models, which take a different approach to approximating PC, outperform PC-SAA. Instead of sampling from the distributions of the random variables, the PW models approximate the ranking-based concordance objective with a win-percentage-based objective which allows us to formulate the problem as a deterministic equivalent. 
PW outperforms in practice because PC-SAA either suffers from (i) too few SAA scenarios causing a loose approximation of PC, or (ii) too many SAA scenarios resulting in  combinatorial explosion and poor convergence (see Appendix~\ref{app:saa_sample_size}). Because PW is represented as a deterministic equivalent, it sidesteps this difficulty as it does not require sampling from the distributions.

\begin{figure}[ht]
    \centering
	\includegraphics[clip,width=\textwidth]{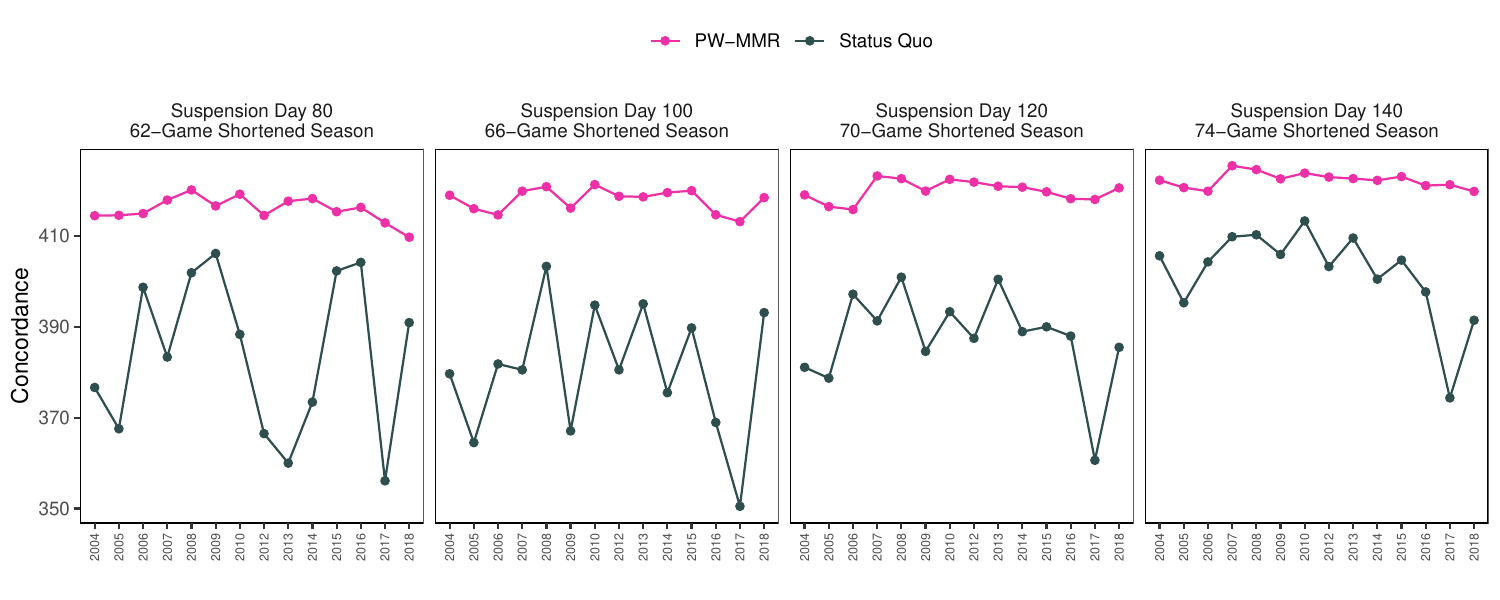}
    \caption{Comparing the best-performing prescriptive model with the suspension day concordance (Status Quo) across 14 NBA seasons.}
    \label{fig:sim_results_status_quo}
\end{figure}

Next, we illustrate the incremental value of playing a shortened season, relative to the ``Status Quo" case of halting the season at the suspension date. Figure \ref{fig:sim_results_status_quo} compares our top-performing prescriptive model (i.e., PW) with ``Status Quo". First, we note that the concordance of ``Status Quo" differs substantially across instances. Indeed, sometimes the ranking on the suspension date is close to the end-of-season ranking (e.g., 2009-10 with suspension day 80 and 2008-09 with suspension day 100), while in other cases the ``Status Quo" ranking is far from the end-of-season ranking (e.g., 2013-14 with suspension day 80 and 2017-18 with suspension days 80 and 100). Second, we observe that playing additional games as chosen by PW significantly improves concordance, regardless of season and suspension date. And finally, as noted previously, as the season is suspended later, all models' rankings converge to the end-of-season ranking, leading to both higher concordance and lower variability in outcomes.

Finally, it is important to note that our methodology is very different from one that takes the ranking at the suspension date and tries to sustain this ranking through the end of the shortened season (had this been the case, the rankings from ``Status Quo" and PW would be similar, and their concordance with the end-of-season ranking as plotted in Figure \ref{fig:sim_results_status_quo} would be much closer). Indeed, although our predictive model uses pre-suspension data to estimate our prescriptive model's parameters, a team with a high win rate pre-suspension will not necessarily have a high win rate post-suspension. The nature of the shortened season and its impact on rankings is complex, owing to the relative difficulty of the schedule before and after suspension (i.e., precisely when a team faces easy or hard-to-beat competitors). Because our prescriptive model aims to produce a shortened season with  ranking similar to that of the full season, if a team had an easy schedule pre-suspension it would generally have a comparatively difficult schedule post-suspension, and our prescriptive model will naturally attempt to maintain this difficulty gradient in the shortened season post-suspension.  
Moreover, it is also interesting to note that the games that PW chooses also tend to be the more ``competitive" games with uncertain outcomes, which are precisely the games that spectators enjoy watching (see Appendix \ref{app:analysis-pw-dqip-model} for details).

Appendix \ref{app:suggestions} provides an illustration of a shortened season produced by our two-phase approach, which is the result of applying our methodology to the suspended 2019-20 NBA regular season.

\subsubsection{Practical implications.}
\label{sec:presc-model-practical-impl}

We now measure the performance of our shortened seasons using several metrics directly tied to a team's ranking. In doing so, we verify that our PW-MMR model produces solutions that not only have high concordance, but are also practically appealing.

As discussed in \cref{subsec:background}, at the end of the regular season, within each conference the top 8 teams make the playoffs and the top 4 receive home court advantage. Moreover, the 5 bottom-ranked teams are given the highest (i.e., double-digit) lottery odds in next year's draft for rookie players. Thus, we define 3 categories of metrics:
(i) ``Playoff teams'' (the top 8 teams in each conference),
(ii) ``Teams with home court advantage'' (the top 4 teams in each conference), and
(iii) ``Teams with double-digit lottery odds'' (the bottom 5 teams overall). Using the same 10,000 simulation replications described in \cref{sec:presc-model-exp-setup}, for each of these 3 metric categories we calculate the proportion of teams in the shortened season that agree with the full season. By \textit{agree}, we mean that we check to see how many teams fall into that category in both the shortened season and the full season. 

For example, imagine that in one particular realization the 5 bottom-ranked teams in the full-season are \{\texttt{Warriors, Cavaliers, Timberwolves, Hawks, Pistons}\}, while the 5 bottom-ranked teams in the shortened-season produced by PW-MMR are \{\texttt{Warriors, Cavaliers, Timberwolves, Hawks, Knicks}\}, and finally the 5 bottom-ranked teams in the shortened-season produced by Greedy are \{\texttt{Warriors, Cavaliers, Timberwolves, Bulls, Knicks}\}.
Here, the PW-MMR model selected $\frac{4}{5} = 80\%$ of the bottom-5 teams correctly while for the Greedy ranking only $\frac{3}{5} = 60\%$ of the bottom-5 teams match those of the full-season. For this realization, we say models PW-MMR and Greedy have 80\% and 60\% \emph{agreement} with the full season, respectively.

Table \ref{table:average_agreement_values} compares the mean agreement percentage of our PW-MMR model with that of Greedy and Status Quo; agreement percentages have been averaged over all 4 suspension dates, 14 seasons, and 10,000 Monte Carlo replications. We can see that if we do not play a shortened season but instead select playoff teams based on the ranking as of the suspension date (Status Quo), then on average 88.96\% of teams that would have made the playoffs if the full season was played do in fact advance to the playoffs; i.e., 88.96\% of playoff teams are chosen correctly. Instead, if we play a shortened season but it is a Greedy rather than optimal one, then on average 94.27\% of the teams that advance to the playoffs are chosen correctly. Finally, if we use PW-MMR to construct a shortened season, on average 95.65\% of the teams advancing to the playoffs would also have been in the playoffs had the full season been played. The other two metrics also indicate significant improvements from using PW-MMR over Greedy and Status Quo. For further details, see 
Appendix \ref{app:practical-implications-boxplots}, which includes boxplots that show the distributions of the three agreement metrics over the 14 NBA seasons, for each of the 4 suspension days.
%






\begin{table}[ht]
\centering
\footnotesize
\renewcommand{\arraystretch}{1.2}
\begin{tabular}{@{}cccc@{}}
\toprule
\textbf{Agreement criteria} & \textbf{PW-MMR} & \textbf{Greedy} & \textbf{Status Quo} \\ \midrule
Playoff teams                          & \textbf{95.65} & 94.27           & 88.96               \\
Teams with home court advantage        & \textbf{92.28} & 89.83           & 79.10               \\
Teams with double-digit lottery odds   & \textbf{91.36} & 89.24           & 78.10               \\ \bottomrule
\end{tabular}
\vspace{0.2cm}
\caption{Mean agreement percentages on simulation for shortened seasons computed using three models.}
\label{table:average_agreement_values}
\end{table}






To further validate the practical performance of our models, we also measured their performance using actual post-suspension game outcomes (i.e., we ran a backtest).
In this experiment, for each season and suspension day, instead of 10,000 simulation replications we have only a single sample path. While this is a good robustness check, the backtest is susceptible to producing noisy outcomes.
We find that, in terms of concordance between the shortened seasons and the full seasons, PW-MMR outperforms Greedy 65\% of the time. Furthermore, when PW-MMR outperforms Greedy, it does so by a larger margin than when Greedy outperforms PW-MMR. Specifically, whenever PW-MMR beats Greedy, its concordance is higher by an average of 3.68 concordant pairs. This is more than twice the number of concordant pairs by which Greedy outperforms PW-MMR (1.57).
Finally, 
Table \ref{table:average_agreement_values_backtest} compares the backtest results for the PW-MMR, Greedy, and Status Quo models using the agreement criteria discussed earlier.

\begin{table}[h]
\centering
\footnotesize
\renewcommand{\arraystretch}{1.2}
\begin{tabular}{@{}cccc@{}}
\toprule
\textbf{Agreement criteria - Backtest} & \textbf{PW-MMR} & \textbf{Greedy} & \textbf{Status Quo} \\ \midrule
Playoff teams                          & \textbf{95.64} & 94.87           & 90.62               \\
Teams with home court advantage        & 90.62 & \textbf{91.07}           & 82.59               \\
Teams with double-digit lottery odds   & \textbf{87.14} & 83.21           & 78.57               \\ \bottomrule
\end{tabular}
\vspace{0.2cm}
\caption{Mean agreement percentages on backtest for shortened seasons computed using three models.}
\label{table:average_agreement_values_backtest}
\end{table}

\subsubsection{Strength of schedule extension.}\label{sec:sos}

A practically appealing extension to our prescriptive model incorporates constraints that additionally ensure a \emph{Strength of Schedule (SoS)} for each team that is not materially reduced from its full-season measure.  There are several mathematical definitions of SoS (see \citealp{NBAStufferSoS}), but essentially SoS is used by each team to quantify the difficulty of its  remaining schedule of games. Pundits on television use SoS for arguments such as, ``Although the Lakers are currently higher-ranked than the Clippers, the Clippers have a higher SoS and so are likely to make up some of this slack and could come out ahead by the end of the season." League managers may wish to assure each team that their SoS is not materially impacted by the shortened season being selected; this motivates the prescriptive model in this subsection. 

As far as we know, there is limited related work on incorporating SoS in a prescriptive model, and the few articles that use such a measure take a dynamic scheduling approach where the schedule is updated according to certain criteria including SoS; see \cite{bouzarth2020dynamically}. As our approach is based on mathematical programming, we choose the \emph{Opponent's Win percentage (OW)} as our SoS metric, as it is both simple to understand and linear in our decision variables.

Following our notation from \cref{sec:prescriptive}, let $\bar{y}_i^0=\frac{y^0_i}{m^0_i}$ denote the win percentage of team $i$ at the time of suspension, where $m^0_i$ is the number of games played by team $i$ pre-suspension. We define team $i$'s strength of schedule in the remainder of the shortened season and full season, respectively, as:
\begin{align}
    OW_i =& \frac{1}{m-m^0_i}\left(\sum\limits_{g\in G^h_i}x_g \bar{y}^0_{j(g)}+\sum\limits_{g\in G^a_i}x_g \bar{y}^0_{i(g)}\right) & \forall i \in T\label{eq:ow_i-def}\\
    \widehat{OW}_i =& \frac{1}{\hat{m}-m^0_i}\left(\sum\limits_{g\in G^h_i}\bar{y}^0_{j(g)}+\sum\limits_{g\in G^a_i} \bar{y}^0_{i(g)}\right)& \forall i \in T \label{eq:owHat_i-def}
\end{align}

The above expressions estimate the average win percentage of all opponents of team $i$ in the remainder of the shortened season and full season, respectively. With slight abuse of notation, $i(g)$ refers to the home team in game $g$, and should not be confused with the focal team $i$. A larger value for $OW_i$ or $\widehat{OW}_i$ implies a more difficult remainder of the season for team $i$, as this means future opponents are harder to beat.

Using the two quantities defined in \eqref{eq:ow_i-def} and \eqref{eq:owHat_i-def}, we modify our PW-DQIP formulation, \eqref{opt:win-percentage-diff-det-obj}--\eqref{opt:win-percentage-diff-det-X}, by adding the constraint \eqref{eq:sos-constraint-PW-SoS} below, which ensures that for each team $i$, the strength of schedule (i.e., average opponents' win percentage) in the remainder of the shortened season is not materially higher than in the remainder of the full season, i.e., within $\epsilon\%$ of the full-season SoS. An extension to our PW-DQIP model incorporating strength-of-schedule constraints is as follows:
\begin{align}
 \text{[PW-SoS]}\quad\min\quad &\sum\nolimits_{i\in T}\left(\left(\mu_i-\hat{\mu}_i\right)^2 +v_i\left(1-\frac{2m}{\hat{m}}\right)+ \hat{v}_i \right) \nonumber \\
    \text{s.t.}\quad & \frac{OW_i - \widehat{OW}_i}{\widehat{OW}_i}
     \, \leq \, \epsilon & \forall i \in T \label{eq:sos-constraint-PW-SoS}\\
    &\eqref{opt:win-percentage-diff-det-mu} \text{--} \eqref{opt:win-percentage-diff-det-X}, \eqref{eq:ow_i-def} \nonumber
\end{align}

Note that in the above, $OW_i$ is an auxiliary decision variable as it depends on the choice of shortened season, while $\widehat{OW}_i$ is a constant since the full season is fixed.

To assess the strength of schedule in our solutions and given we are only interested in cases where the strength of schedule in the shortened season is larger than that of full season, we introduce the metric ``Strength of Schedule Discrepancy (SSD)'' which takes the \emph{positive difference} between strength of schedule in the shortened and full seasons, defined as:
\begin{align}
    \text{SSD} = \frac{1}{n} \sum_{i \in T} \max \left\{\frac{OW_i - \widehat{OW}_i}{\widehat{OW}_i}, 0\right\}.\label{eq:pos-diff-def}
\end{align}

\begin{figure}[t]
    \centering
    \includegraphics[clip,width=\textwidth]{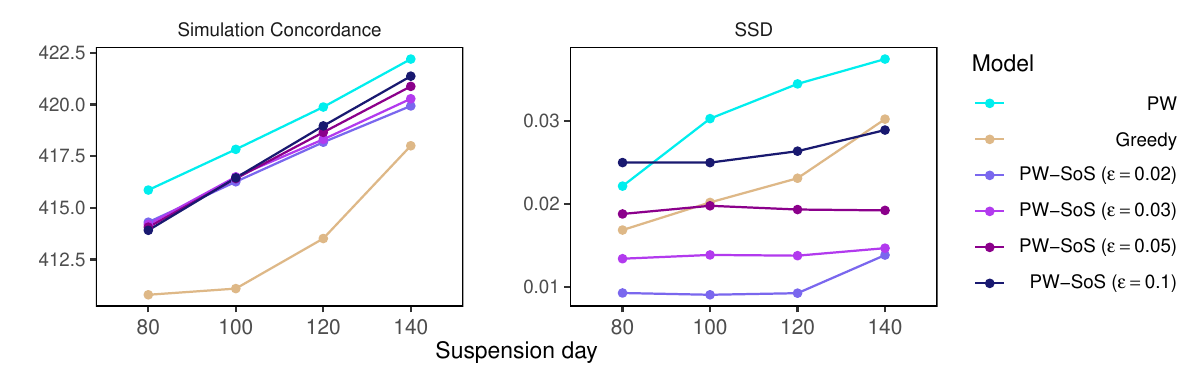}
    \caption{Comparison of PW-SoS (with 4 different values for $\epsilon$) to PW and Greedy based on simulation concordance (left panel), and SSD (right panel).}
    \label{fig:OWP_overall_ssd_sim}
\end{figure}

Figure \ref{fig:OWP_overall_ssd_sim} compares models PW (without SoS constraints), Greedy and PW-SoS with 4 choices for $\epsilon$ ranging from 2\% to 10\%. The left panel plots concordance as measured in our Monte Carlo simulation over 10,000 replications. Here, we see that the addition of SoS constraints sacrifices some amount of concordance, but in general PW-SoS still performs significantly better than Greedy. In the right panel, we investigate the sensitivity to the $\epsilon$ parameter, and show that for $\epsilon$ below 3\% we can produce solutions using PW-SoS that have better (lower) SSD than both PW and Greedy.
Each dot in Figure \ref{fig:OWP_overall_ssd_sim} represents an average value across 14 NBA seasons. Appendix \ref{app:SoS-lineplots} includes the corresponding plots for individual NBA seasons.




\section{Conclusions and Future Research} \label{sec:conclusion}

Professional sports leagues may be suspended due to various reasons, requiring the league to select which games to play in a shortened season. In this paper, we proposed a two-phase analytics approach for this problem. 
In phase one, we predicted game outcomes using a composite binary classifier, with a particular functional form for each season and for each suspension day, chosen based on LogLoss values. In phase two, we used stochastic optimization techniques to prescribe a data-driven decision which maximizes the expected similarity between the ranking at the end of the shortened season and the full season had it been played in full. 

To solve one of our stochastic optimization problems (PW), we proposed three solution methodologies: i) a deterministic equivalent reformulation (i.e., PW-DQIP), ii) a Frank-Wolfe decomposition algorithm (i.e., PW-FW) which significantly reduces the running time of  PW-DQIP without sacrificing solution quality, and iii) a robust reformulation of  PW-DQIP designed to handle misspecification in the input data (i.e., PW-MMR).
For our second model (PC), we proposed approximation schemes (MVP and SAA), as well as variable fixing techniques. Our PC-SAA model approximates the distributions in PC but has an exact objective, while our PW models use the exact distribution but approximate PC's objective.

We evaluated our models' solutions using Monte Carlo simulation.
Our computational experiments show PW outperforms PC-SAA even for reasonably-large 50-scenario instances. This suggests that PW (and specifically the PW-MMR variant) is the recommended prescriptive model to maximize concordance. Finally, we verified that the higher-concordance solutions provided by PW-MMR outperform our benchmarks, leading to a higher agreement between shortened season and counterfactual full season in terms of the number of teams that (a) make the playoffs, (b) receive home court advantage, and (c) receive double-digit rookie draft lottery odds. We ran a backtest as a robustness check, and also provided a model extension (PW-SoS) that ensures each team's strength-of-schedule is not materially impacted by our choice of shortened season.

We envision several directions for future research.
First, one potential improvement could come from considering more sophisticated loss functions for the predictive models that are custom-tailored to the specific needs of the downstream prescriptive models.
Second, apart from concordance, other considerations (e.g., generated revenue, travel cost and distances, broadcasting restrictions, venue availability) may also be relevant, suggesting an alternative multi-criteria decision-making approach. Third,  when faced with multiple stoppages during a single season, selecting the optimal subset of games between any two stoppages can be modeled as a dynamic stochastic optimization problem with a learning component. Finally, large-scale stochastic optimization techniques (e.g., progressive hedging) may be designed to tackle SAA with a larger sample.




\bibliographystyle{informs2014} 
\bibliography{ref.bib}

\normalsize

%
%
%

\newpage
\begin{APPENDICES}

\section{Solution Methods for the Stochastic Model PC}\label{app:approx-solution-methods-pm-PC}
In this section, we introduce two solution methodologies based on model PC, introduced in \cref{sec:model-pc}. Before elaborating on the solution methodologies, we remark that our formulation PC may be viewed as a stochastic program with recourse, where the $x$-variables are first-stage variables (for which there is only one choice to be made) and the $\vwp$- and $z$-variables are second-stage ``recourse" variables (for which there is one such variable for each possible outcome $\xi$).  Note, however, that in our application there is no true recourse.  Rather, $\vwp$ and $z$ are auxiliary variables whose purpose is to linearize the objective function. In the next section, as the full stochastic optimization problems are too large to solve directly, we introduce two methods which approximately solve PC.

\subsection{Mean Value Approximation}\label{subsec-mvp}
Replacing all random parameters in a stochastic optimization problem by their expected values yields a deterministic problem known as the Mean Value Problem (MVP). In our case, we may produce MVP for PC  by replacing the random variables $W_g$ which represent the outcome of each game $g$ with their means $p_g=\mathbb{E}[W_g]$. The $\vwp$-variables are then interpreted as expected values over all outcomes $\xi \in \Xi$, given the shortened season $\xx$, and $z$ variables capture relative positions of teams according to $\vwp$. The MVP corresponding to PC is:

\begin{align}
    \mbox{[PC-MVP]}\;\max\quad & \;\sum_{i\in T}\sum_{j\in T:j>i}\left( z_{ij}\hat{z}_{ij}+(1- z_{ij})(1-\hat{z}_{ij})\right)\label{opt:mvp-kentall-tau-obj} \\
    \mbox{s.t.}\quad&\vwp_i=\frac{1}{m}\left(\vwp^0_i + \sum_{g\in G^h_i}p_g x_g+ \sum_{g\in G^a_i}(1-p_g) x_g\right) & \forall i\in T \label{opt:mvp-kentall-tau-yi-def}\\ 
    &z_{ij}\ge \vwp_i - \vwp_j\ge z_{ij}-1 & \forall i,j\in T: i<j \label{opt:mvp-kentall-tau-zij-def}\\
    &z_{ij}\in \{0,1\} & \forall i,j\in T: i<j\label{opt:mvp-kentall-tau-zij-domain}\\
    & \xx\in X.\label{opt:mvp-kentall-tau-x-domain}
\end{align}

\subsection{Sample Average Approximation}
\label{subsec-SAA}
Sample Average Approximation (SAA) is a Monte Carlo simulation-based technique for approximating stochastic optimization problems \citep{kleywegt2002sample},
Let $\mathcal{S}=\{\xi^{(1)},\xi^{(2)},\dots,\xi^{(|\mathcal{S}|)}\}$ be an independently and identically distributed random sample of $\xi$.
SAA reduces the size of the problem by approximating the expected value in the objective function with the sample average function. We use the superscript $s$ to reference the second-stage variables and random parameters under scenario $s\in \mathcal{S}$. For instance, under scenario $s$, $W^{(s)}_g$ refers to the outcome of game $g$, $\hat{\vwp}^{(s)}_i$ refers to the win percentage of team $i$ at the end of the full season, and $\vwp^{(s)}_i$ refers to the decision variable for the win percentage of team $i$ at the end of the shortened season.
We construct the SAA counterpart of the stochastic program PC by replacing the full set of outcomes $\Xi$ with the sample set $\mathcal{S}$. 
\begin{align}
    \mbox{[PC-SAA]}\;\max\quad &\frac{1}{|\mathcal{S}|}\sum_{s\in \mathcal{S}} \sum_{i\in T}\sum_{j\in T:j>i}\left( z^{(s)}_{ij}\hat{z}^{(s)}_{ij}+(1- z^{(s)}_{ij})(1-\hat{z}^{(s)}_{ij})\right)\label{opt:saa-kentall-tau-obj} \\
    \mbox{s.t.}\quad & \vwp^{(s)}_i=\frac{1}{m}\left(\vwp^0_i + \sum_{g\in G^h_i}W^{(s)}_g x_g+ \sum_{g\in G^a_i}(1-W^{(s)}_g) x_g\right) & \forall i\in T,\forall s\in \mathcal{S} \label{opt:saa-kentall-tau-yi-def}\\
    &z^{(s)}_{ij}\ge \vwp^{(s)}_i - \vwp^{(s)}_j\ge z^{(s)}_{ij}-1 & \forall i,j\in T:i<j,\forall s\in \mathcal{S} \label{opt:saa-kentall-tau-zij-def}\\
    &z^{(s)}_{ij}\in \{0,1\} & \forall i,j\in T:i<j,\forall s\in \mathcal{S}\label{opt:saa-kentall-tau-zij-domain}\\
    & \xx\in X.\label{opt:saa-kentall-tau-x-domain}
\end{align}

As the sample size increases, the optimal solution and the optimal value of the SAA problems converge to their `true' stochastic counterparts with probability one \citep{kleywegt2002sample}. 

\subsection{Variable Fixing and Preprocessing}\label{app:variable-fixing}
We may improve the computational efficiency of both the SAA and MVP counterparts of PC by fixing certain variables at their optimal values and eliminating redundant constraints, as described by the following proposition. 

\begin{proposition}\label{prop:variable-fixing}
Let $\tilde\xi$ be an arbitrary realization or expected value of $\xi$. For each team $i$, sort $G^h_i$ and $G^a_i$ in non-decreasing order of $W(\tilde\xi)$. Let $U^h_i$ and $L^h_i$ be the summation of $W_g(\tilde\xi)$ values corresponding to the first and last $m^h_i$ games in $G^h_i$, respectively. Similarly, let $U^a_i$ and $L^a_i$ be the summation of $W_g(\tilde\xi)$ values corresponding to the first and last $m^a_i$ games in $G^a_i$, respectively. Define $\vwp^U_i=\frac{1}{m}(\vwp^0_i+U^h_i+m^a_i-L^a_i)$ and $\vwp^L_i=\frac{1}{m}(\vwp^0_i+L^h_i+m^a_i-U^a_i)$ to be the optimistic and pessimistic win percentages of team $i$ under $\tilde\xi$, respectively. For each pair of teams $(i,j)$:
\begin{enumerate}[(i)]
    \item If $\vwp^L_i-\vwp^U_j> 0$, then $z_{ij}(\tilde\xi)=1$, and the corresponding linking constraints are redundant.
    \item If $\vwp^U_i-\vwp^L_j< 0$, then $z_{ij}(\tilde\xi)=0$, and the corresponding linking constraints are redundant.
\end{enumerate}
\end{proposition}
\proof{Proof.} 
Using the definition of the $z$-variables, the statements follow from $z_{ij}(\tilde\xi)\ge \vwp_i(\tilde\xi)-\vwp_j(\tilde\xi)\ge \vwp^L_i-\vwp^U_j>0$, and $ 0>\vwp^U_i-\vwp^L_j\ge \vwp_i(\tilde\xi)-\vwp_j(\tilde\xi)\ge z_{ij}(\tilde\xi)-1$, respectively.
\Halmos
\endproof

Table \ref{table:elimination} summarizes the results of applying our variable fixing technique, introduced in Proposition \ref{prop:variable-fixing}, in PC variants.
The high percentages under the columns ``Percentage'' in Table \ref{table:elimination} highlight the effectiveness of the variable fixing technique in eliminating a large proportion of the $z$-variables across different scenarios in both MVP and SAA. More importantly, the technique is able to eliminate between 8,000--17,000 binary variables in the SAA problems with 50 scenarios. 
We also observe that as the suspension day increases (i.e., the season is suspended later), more pairs of teams become impossible to switch ranking positions, given the limited number of remaining games. For instance, when the season is suspended on day 140, PC-SAA has control over only 20\% of these binary variables in a 74-game shortened season, with the remaining 80\% of the variables fixed (i.e., eliminated).

\begin{table}[t!]
\centering
\footnotesize
\begin{tabular}{@{}lllll@{}}
\toprule
\multicolumn{1}{c}{\textbf{Sus. Day}}  & \multicolumn{2}{c}{MVP}                                        & \multicolumn{2}{c}{SAA}                                        \\ \cmidrule(l){2-5} 
\multicolumn{1}{c}{(GT)} & \multicolumn{1}{c}{Percentage} & \multicolumn{1}{c}{Variables} & \multicolumn{1}{c}{Percentage} & \multicolumn{1}{c}{Variables} \\ \midrule
80 (62)                     & 62.9\%                        & 273.7                         & 38.6\%                        & 8,393.0                        \\
100 (66)                     & 70.1\%                        & 304.9                         & 53.9\%                        & 11,722.8                        \\
120 (70)                     & 79.6\%                        & 346.4                         & 65.5\%                        & 14,240.7                          \\
140 (74)                     & 88.1\%                        & 383.2                         & 80.1\%                        & 17,423.0                          \\ \midrule
Average & 75.2\%                        & 327.1                         & 59.5\%                        & 12,944.9 \\
\bottomrule
\end{tabular}
\caption{Number of $z$-variables eliminated by our variable fixing technique, reported for different suspension days (80, 100, 120, 140) and target number of games/team (62,66,70,74).  Results are averaged over 14 seasons.}
\label{table:elimination}
\end{table}

\subsection{Tuning the Sample Size for SAA}\label{app:saa_sample_size}
Here we analyze the impact of the number of scenarios on our PC-SAA model.
Figure~\ref{fig:prescripive_saa} presents the performance of PC-SAA across six choices of sample size $|\mathcal{S}|\in\{5,10,25,50,100,200\}$.
Each boxplot corresponds to 14 values for 14 NBA seasons, assuming a suspension day 100 and 66 as the target number of games in the shortened season.
The concordance values are obtained after evaluating the solution $\xx$ proposed by each model on 1,000 randomly-generated scenarios.
\begin{figure}[ht]
    \centering
    \includegraphics[clip,width=1\textwidth]{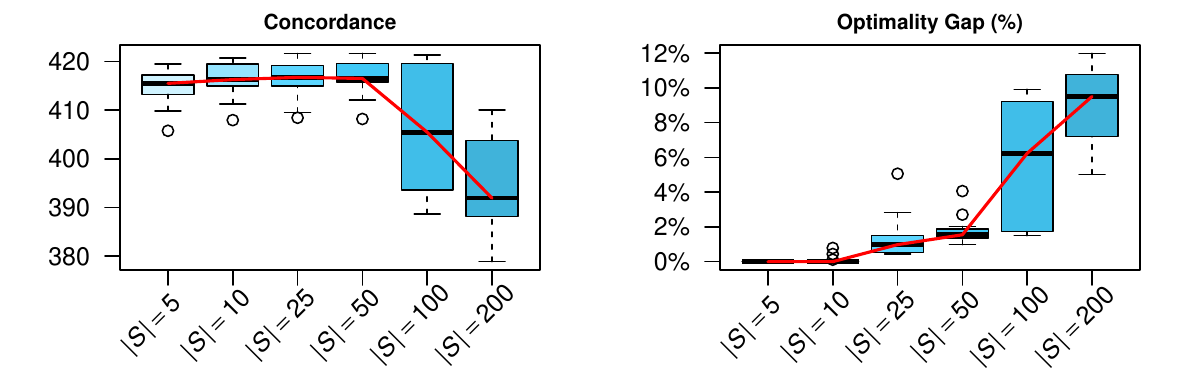}
    \caption{Performance of the SAA algorithm across different choices of sample size.}
    \label{fig:prescripive_saa}
\end{figure}
The panel on the right presents the optimality gaps of the SAA problems after reaching a time limit of 3600 seconds. As the number of scenarios increases, one should expect to obtain a closer approximation of the true stochastic problem via SAA. However, a larger sample amounts to solving a more challenging SAA problem. As depicted in Figure~\ref{fig:prescripive_saa}, initially as the sample size increases, the quality of the solution improves in the simulation phase. However, after surpassing 50 scenarios, the SAA becomes computationally intractable, leading to larger optimality gaps and a degradation in the quality of the solution. Hence, the trade-off between the quality of the solution and the runtime is balanced at 50 scenarios. Thus, we select 50 scenarios for the SAA counterpart of our PC model in our main experiments.

\section{Proof of Theorems and Propositions}\label{app:proofs}

\proof{Proof of Proposition~\ref{prop:kendall_euclidean}.} Let us define the Manhattan distance between rankings $r$ and $\hat{r}$ as 
\begin{align*}
    d_{\text{M}}(r,\hat{r}) = \sum\nolimits_{i\in T} |r_i-\hat{r}_i|.
\end{align*}
We establish the relationship between $d_{\text{E}}(r,\hat{r})$ and  $\tau_{\text{C}}(r,\hat{r})$ described in \eqref{eq:kendall_euclidean} by showing that
\begin{align}
    \sqrt{2 \,d_{\text{E}}(r,\hat{r})} \le & d_{\text{M}}(r,\hat{r}), \text{ and} \label{eq:manhattan_eucliean}\\
    d_{\text{M}}(r,\hat{r})\le & n(n-1) - 2 \tau_{\text{C}}(r,\hat{r}). \label{eq:manhattan_kendall}
\end{align}
To prove inequality \eqref{eq:manhattan_eucliean}, define $d_i=r_i-\hat{r}_i$, and note that $d_{\text{E}}(r,\hat{r})=\sum_{i\in T}d_i^2$, and $d_{\text{M}}(r,\hat{r})=\sum_{i\in T}|d_i|$. Observe that $\sum_{i\in T}d_i=\sum_{i\in T}r_i-\sum_{i\in T}\hat{r}_i=0$, which, using triangle inequality, implies that for each $i\in T$:
\begin{align*}
    |d_i|=|-\sum\nolimits_{j\ne i} d_j|\le \sum\nolimits_{j\ne i} |d_j| = d_{\text{M}}(r,\hat{r})- |d_i| \Rightarrow d_{\text{M}}(r,\hat{r}) \ge 2|d_i|.
\end{align*}
Consequently, expanding $d^2_{\text{M}}(r,\hat{r})$, we obtain
\begin{align*}
    d^2_{\text{M}}(r,\hat{r})=d_{\text{M}}(r,\hat{r})\sum\nolimits_{i\in T} |d_i|\ge 2\sum\nolimits_{i\in T} |d_i|^2 = 2\sum\nolimits_{i\in T} d_i^2 = 2d_{\text{E}}(r,\hat{r}),
\end{align*}
which establishes \eqref{eq:manhattan_eucliean}.
To prove \eqref{eq:manhattan_kendall}, let us define $z_{ij}=\mathbb{I}(r_i<r_j)$ and $\hat{z}_{ij}=\mathbb{I}(\hat{r}_i<\hat{r}_j)$ for each pair $(i,j)$, where $\mathbb{I}(\cdot)$ is the indicator function. Observe that $1-|z_{ij}-\hat{z}_{ij}|$ is 1 when pair $(i,j)$ is concordant across $r$ and $\hat{r}$ and is 0 otherwise. Hence, we can compute $\tau_{\text{C}}(r,\hat{r})$ using
\begin{align}
    \tau_{\text{C}}(r,\hat{r}) = \frac{1}{2} \sum\nolimits_{i\in T}\sum\nolimits_{j\in T:j\ne i}(1-|z_{ij}-\hat{z}_{ij}|)=\frac{1}{2}n(n-1)-\frac{1}{2} \sum\nolimits_{i\in T}\sum\nolimits_{j\in T:j\ne i}|z_{ij}-\hat{z}_{ij}|. \label{eq:kendall_z}
\end{align}
On the other hand, $r_i=n-\sum_{j\ne i}z_{ij}$, since $\sum_{j\ne i}z_{ij}$ counts how many teams have a worse rank than $i$. Similarly, $\hat{r}_i=n-\sum_{j\ne i}\hat{z}_{ij}$, yielding $r_i-\hat{r}_i=\sum_{j\in T:j\ne i}(z_{ij}-\hat{z}_{ij})$. Hence, \eqref{eq:kendall_z} and triangle inequality imply
\begin{align*}
    d_{\text{M}}(r,\hat{r})=\sum\nolimits_{i\in T}|r_i-\hat{r}_i| = \sum\nolimits_{i\in T}\left|\sum\nolimits_{j\in T:j\ne i}(z_{ij}-\hat{z}_{ij})\right|\le \sum\nolimits_{i\in T}\sum\nolimits_{j\in T:j\ne i}|z_{ij}-\hat{z}_{ij}|=n(n-1) - 2 \tau_{\text{C}}(r,\hat{r}),
\end{align*}
which establishes \eqref{eq:manhattan_kendall} and completes the proof. \Halmos\endproof

~

\begin{lemma}\label{lemma:max-l2-permutation}
    Maximum quadratic Euclidean distance between any two permutations of $\{1,\dots,n\}$ is $\frac{n}{3}(n^2-1)$.
\end{lemma}
\proof{Proof.}
Let $P$ denote the set of all permutations of $N\coloneqq \{1,\dots,n\}$. Our goal is to find permutations $p\in P$ and $q\in P$ which maximize $\sum\nolimits_{i\in N}^n(p_i-q_i)^2$. Note that, without loss of generality, we can fix $p=(1,2,\dots,n)$, and restate the problem as
\begin{align}
    \max_{q\in P}\; \sum\nolimits_{i\in N} (q_i-i)^2. \label{eq:max-permutation}
\end{align}
We first note that setting $q=(n,n-1,\dots,1)$ (i.e., the exact reverse of $p$) yields $\sum\nolimits_{i\in N}(q_i-i)^2 = \sum\nolimits_{i\in N} \left((n+1-i)- i\right)^2 = \sum\nolimits_{i\in N} \left(2i - n - 1\right)^2 = 4 \sum\nolimits_{i\in N} i^2 \, - \, 4(n+1) \sum\nolimits_{i\in N} i+n(n+1)^2= \frac{n}{3}(n^2-1)$. 

Next we use LP duality to show that this lower bound is tight. 
Note that we may state \eqref{eq:max-permutation} as the following assignment problem
\begin{align*}
    \max\quad & \sum\nolimits_{i\in N}\sum\nolimits_{j\in N} x_{i,j} (i-j)^2\\
    \text{s.t.}\quad & \sum\nolimits_{j\in N} x_{i,j}=1 &&\forall i\in N\\
    & \sum\nolimits_{i\in N} x_{i,j}=1 &&\forall j\in N\\
    & x_{i,j}\ge 0 &&\forall i,j\in N,
\end{align*}
which can be stated in the dual form as
\begin{align*}
    \min\quad & \sum\nolimits_{i\in N} \alpha_i + \sum\nolimits_{j\in N} \beta_j\\
    \text{s.t.} \quad & \alpha_i + \beta_j \ge (i-j)^2 \qquad \forall i,j\in N.
\end{align*}
It is not difficult to verify that setting 
\begin{align*}
    \alpha_i=\beta_i = \begin{cases}
        \frac{1}{2}(n-i)^2 & \text{ if } i \le \frac{n}{2}\\
        \frac{1}{2}(i-1)^2 & \text{ if } i> \frac{n}{2}
    \end{cases} \qquad \forall i\in N
\end{align*}
satisfies $\alpha_i + \beta_j = \alpha_i + \alpha_j \ge (i-j)^2$ for each $i$ and $j$, and yields $\sum_{i\in N} \alpha_i + \sum_{j\in N} \beta_j = 2\sum_{i\in N} \alpha_i = \frac{n}{3}(n^2-1)$. Hence, by strong duality, the optimal value for \eqref{eq:max-permutation} is $\frac{n}{3}(n^2-1)$.
\Halmos\endproof

~

\proof{Proof of Proposition \ref{prop:euclidean_rank_winpercentage}.} The statement holds when win percentages are identical, which results in 0 on both sides. Now, assuming that win percentages are not identical, there exists team $j$ such that $\vwp_j(\xx,\xi)\ne \hat{\vwp}_j(\xi)$. 
Note that $\hat{\vwp}_j(\xi)\in \{0,\frac{1}{\hat{m}},\frac{2}{\hat{m}},\dots,1\}\subseteq \{0,\frac{1}{L},\frac{2}{L},\dots,1\}$, since $\hat{\vwp}_j(\xi)$ is the number of wins for team $j$ in the full season divided by $\hat{m}$. Similarly, $\vwp_j(\xx,\xi)\in \{0,\frac{1}{m},\frac{2}{m},\dots,1\}\subseteq \{0,\frac{1}{L},\frac{2}{L},\dots,1\}$.
Therefore, the closest that $\vwp_j(\xx,\xi)$ and $\hat{\vwp}_j(\xi)$ can get and still be different is $\frac{1}{L}$; thus for non-identical win percentage vectors we have
\begin{align}
    \frac{1}{L^2}\le \sum\nolimits_{i\in T}\left(\vwp_i(\xx,\xi)-\hat{\vwp}_i(\xi)\right)^2 \label{eq:min_winpercentage}
\end{align}
On the other hand, by Lemma \ref{lemma:max-l2-permutation} we have
\begin{align}
    d_{\text{E}}(r(\xx,\xi),\hat{r}(\xi)) \le \frac{n}{3}(n^2-1). \label{eq:max_euclidean}
\end{align}

Multiplying both sides of \eqref{eq:min_winpercentage} and \eqref{eq:max_euclidean} and rearranging the resulting inequality yields
\begin{align*}
    d_{\text{E}}(r(\xx,\xi),\hat{r}(\xi)) \le \frac{n}{3}(n^2-1) L^2 \sum\nolimits_{i\in T}\left(\vwp_i(\xx,\xi)-\hat{\vwp}_i(\xi)\right)^2,
\end{align*}
which shows existence of $D\le \frac{n(n^2-1)}{3} L^2$. \Halmos
\endproof

~

\proof{Proof of Theorem \ref{theorem-pw}.} 
Using identity $\mathbb{E}[X^2]=\mathbb{E}[X]^2+\mathbb{V}[X]$, with $\mathbb{V}(\cdot)$ denoting variance, we obtain
\begin{align*}
    \mathbb{E}_{\xi}\left[\sum\nolimits_{i\in T}\left(\vwp_i(\xi)-\hat{\vwp}_i(\xi)\right)^2\right]=\sum_{i\in T}\mathbb{E}_{\xi}\left[\left(\vwp_i(\xi)-\hat{\vwp}_i(\xi)\right)^2\right]=\sum_{i\in T}\mathbb{E}_{\xi}\left[\vwp_i(\xi)-\hat{\vwp}_i(\xi)\right]^2+\sum_{i\in T}\mathbb{V}_{\xi}\left[\vwp_i(\xi)-\hat{\vwp}_i(\xi)\right].
\end{align*}
Clearly, $\mathbb{E}_{\xi}\left[\vwp_i(\xi)-\hat{\vwp}_i(\xi)\right]=\mu_i-\hat{\mu}_i$. Moreover, given the definition of $\vwp_i(\xi)$ and $\hat{\vwp}_i(\xi)$, we have
\begin{align}
    \vwp_i(\xi)-\hat{\vwp}_i(\xi)=&(\frac{1}{m}-\frac{1}{\hat{m}})\vwp^0_i +\sum\nolimits_{g\in G_i^h} W_g(\xi)\left(\frac{1}{m}x_g-\frac{1}{\hat{m}}\right)+\sum\nolimits_{g\in G_i^a} (1-W_g(\xi))\left(\frac{1}{m}x_g-\frac{1}{\hat{m}}\right)\nonumber \\
   \Rightarrow \mathbb{V}_{\xi}\left[\vwp_i(\xi)-\hat{\vwp}_i(\xi)\right] =&\sum\nolimits_{g\in G_i^h\cup G_i^a}p_g(1-p_g)\left(\frac{1}{m}x_g-\frac{1}{\hat{m}}\right)^2, \label{eq:pw-dqip-proof-1}
\end{align}
where we have used $\mathbb{V}_{\xi}\left[W_g(\xi)\right]=\mathbb{V}_{\xi}\left[1-W_g(\xi)\right]=p_g(1-p_g)$. Given that $x_g\in\{0,1\}$, we have
\begin{align}
    \left(\frac{1}{m}x_g-\frac{1}{\hat{m}}\right)^2 = \frac{1}{m^2}x_g^2-\frac{2}{m\hat{m}}x_g+\frac{1}{\hat{m}^2} = \frac{x_g}{m^2}\left(1-\frac{2m}{\hat{m}}\right) +\frac{1}{\hat{m}^2},\label{eq:pw-dqip-proof-2}
\end{align}
where we have used $x_g^2=x_g$.
Replacing \eqref{eq:pw-dqip-proof-2} into \eqref{eq:pw-dqip-proof-1} yields
\begin{align*}
   \mathbb{V}_{\xi}\left[\vwp_i(\xi)-\hat{\vwp}_i(\xi)\right] =&\sum\nolimits_{g\in G_i^h\cup G_i^a}p_g(1-p_g)\left(\frac{x_g}{m^2}\left(1-\frac{2m}{\hat{m}}\right) +\frac{1}{\hat{m}^2}\right)\\
   =& \left(1-\frac{2m}{\hat{m}}\right)\frac{1}{m^2}\sum_{g\in G_i^h\cup G_i^a}p_g(1-p_g)x_g + \frac{1}{\hat{m}^2}\sum_{g\in G_i^h\cup G_i^a}p_g(1-p_g)=\left(1-\frac{2m}{\hat{m}}\right)v_i+\hat{v}_i,
\end{align*}
which completes the proof. \Halmos
\endproof

~

\proof{Proof of Proposition \ref{prop:tum}.} 
As illustrated in Figure~\ref{fig:multigraph}, the set $X$ corresponds to a bipartite multigraph with $2n$ nodes ($n$ home teams and $n$ away teams) and each game $g\in G$ corresponds to an edge of unit capacity between home team $i(g)$ and away team $j(g)$. The coefficient matrix of $X$ is the incidence matrix of this bipartite multigraph, which is totally unimodular \citep[see e.g.,][]{yannakakis1985class}. \Halmos
\endproof

\begin{figure}[ht]
    \centering
    \includegraphics[width=0.6\textwidth]{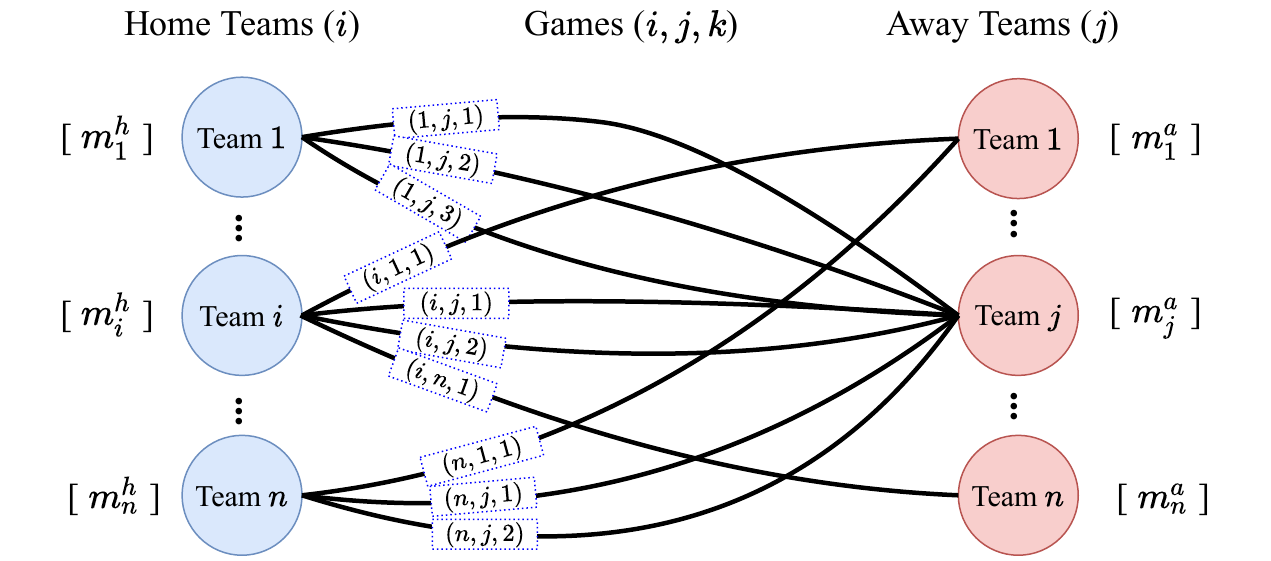}
    \caption{Set of feasible schedules $X$ corresponds to a bipartite multigraph.}
    \label{fig:multigraph}
\end{figure}

\section{Analysis of the PW-DQIP Model}\label{app:analysis-pw-dqip-model}

Let us take a closer look at the mathematical formulation of the PW-DQIP model, and its objective function in particular. As discussed in Section 5, PW-DQIP approximates the objective function in equation (8) while formulating a deterministic quadratic integer program. Once, we expand the objective function in PW-DQIP, (12), the resulting expression consists of i) \emph{L2 Norm} of the difference between the vector of win percentages in the shortened season and a target value (i.e., win percentages in the full season), ii) variance of the win percentages in the shortened season multiplied by a constant factor ($\alpha$), and iii) variance of the win percentages in the full season. Note that the last term $\hat{v}_i$ is constant, as it does not include any decision variables.
{\footnotesize
\begin{align}
    \mathbb{E}_{\xi}\left[\sum\nolimits_{i\in T}\left(\vwp_i(\xi)-\hat{\vwp}_i(\xi)\right)^2\right] &= \sum\nolimits_{i\in T}\left(\left(\mu_i-\hat{\mu}_i\right)^2 +v_i\left(1-\frac{2m}{\hat{m}}\right)+ \hat{v}_i \right)\label{eq:pw-dqip-obj-original}\\
    &= \sum\nolimits_{i\in T}\left(\overbrace{\left(\mu_i-\hat{\mu}_i\right)^2}^{\text{L2 Norm}}+\overbrace{\left(1-\frac{2m}{\hat{m}}\right)}^{\text{Constant $\alpha$}} \,\overbrace{v_i}^{\text{Variance}}+ \overbrace{\hat{v}_i}^{\text{Constant}} \right) \label{eq:pw-dqip-obj-decomposed}
\end{align}}
Omitting the third (constant $\hat{v}_i$) term in \eqref{eq:pw-dqip-obj-decomposed}, the model PW-DQIP equivalently optimizes the following function:
\begin{align}
   \sum\nolimits_{i\in T}\left(\left(\mu_i-\hat{\mu}_i\right)^2 \,+ \, \alpha \, v_i\right), \label{eq:pw-dqip-obj-no-constant}
\end{align}
where:
\begin{align}
    \alpha = 1-\frac{2m}{\hat{m}}.\label{eq:alpha-m-mHat-def}
\end{align}
Note that the number of games per team in the full regular season is fixed at 82 per the structure of the NBA (i.e., $\hat{m} = 82$). Depending on the target number of games in the shortened season ($m$), parameter $\alpha$ could be positive, negative, or zero. According to \eqref{eq:alpha-m-mHat-def}, parameter $\alpha$ is zero when $m = \frac{\hat{m}}{2} = \frac{82}{2} = 41$. In other words, when the target number of games in the shortened season is exactly half of the length of the full season, the coefficient of the variance term in the objective function disappears. The variance term has a positive coefficient when $m < \frac{\hat{m}}{2}$ (shortened season is ``short''), and a negative coefficient when $m > \frac{\hat{m}}{2}$ (shortened season is ``long''). Given the shortened season is long ($m > \frac{\hat{m}}{2}$) in all of the suspension instances in our experiments ($m \in \{62, 66, 70, 74\}$), the PW-DQIP model favors games with higher variance, as the objective function in \eqref{eq:pw-dqip-obj-no-constant} is to be minimized and the coefficient of the variance term ($\alpha$) is negative.

Now, let us revisit our argument focusing on individual games and how the inclusion or exclusion of a game contributes to the variance term in the objective function. The objective function in PW-DQIP can be reformulated as:
\begin{align*}
    &\sum_{i\in T}\left(\left(\mu_i-\hat{\mu}_i\right)^2 +v_i\left(1-\frac{2m}{\hat{m}}\right)+ \hat{v}_i \right)=\sum_{i\in T}\hat{v}_i+\sum_{i\in T}\left(\mu_i-\hat{\mu}_i\right)^2 + 2\left(1-\frac{2m}{\hat{m}}\right)\frac{1}{m^2}\sum_{g\in G}p_g(1-p_g)x_g
\end{align*}
For a game $g \in G$, the variance for the predicted probability $p_g$, following the Bernoulli distribution, is $p_g (1 - p_g)$. It is not difficult to see that the variance term is maximized when our estimated $p_g$ is closer to 0.5. We can also define a \emph{sharpness} metric as $\max \{p_g, 1- p_g\}$ to assess how sharp the estimated probabilities are. A higher variance in game outcomes
leads to lower sharpness and vice versa. In other words, higher variance in the predicted probabilities will favor more evenly matched games, while a lower variance will favor more one-sided matches.
Table \ref{table:short-long-shortened-season} summarizes our discussion in this section.

\begin{table}[t!]
\centering
\footnotesize
\renewcommand{\arraystretch}{1.2}
\scalebox{0.9}{
\begin{tabular}{cccC{6cm}}
\toprule
Target \# games ($m$)         & Shortened season & Variance coefficient ($\alpha$) & Desired games by the model PW-DQIP                                             \\ \midrule
$m < \frac{\hat{m}}{2}$ & Short            & Positive                       & One-sided matchups     \\
$\bm{m > \frac{\hat{m}}{2}}$ & \textbf{Long}             & \textbf{Negative}                       & \textbf{Evenly matched games} \\ \bottomrule
\end{tabular}}
\vspace{0.2cm}
\caption{Comparing short vs. long shortened seasons and the practical implications in the model PW-DQIP.}
\label{table:short-long-shortened-season}
\end{table}

To see the effect of suspension day on the average variance of the outcomes of the selected games compared to the average variance of the excluded games, we further note that:
\begin{align*} 
    &\sum_{i\in T}\hat{v}_i+\sum_{i\in T}\left(\mu_i-\hat{\mu}_i\right)^2 + 2\left(1-\frac{2m}{\hat{m}}\right)\frac{1}{m^2}\sum_{g\in G}p_g(1-p_g)x_g=\sum_{i\in T}\hat{v}_i+\sum_{i\in T}\left(\mu_i-\hat{\mu}_i\right)^2 + 2\left(1-\frac{2m}{\hat{m}}\right)\frac{G_1}{m^2}\bar{v}
\end{align*}
where $G_1$ is the number of post-suspension games included in the shortened season (note that $G_1= \frac{1}{2}\sum_{i\in T} (m_i^a+m_i^h)=\sum_{g\in G}x_g$ for any feasible shortened season $\xx\in X$), and $\bar{v}$ is the average variance of these games, i.e.:
\begin{align*}
    \bar{v} = \frac{1}{G_1}\sum_{g\in G}p_g(1-p_g)x_g.
\end{align*}
Thus the coefficient $2\left(1-\frac{2m}{\hat{m}}\right)\frac{G_1}{m^2}$ governs the trade-off between the contribution of the average variance of the selected games and the L2 term to the objective function.
As depicted in Figure \ref{fig:pd-dqip-var-coeffs}, this coefficient increases (its magnitude decreases) as suspension day increases. As a result, relative to the L2-norm component of the objective function, the variance term has a lower weight, and one should expect lower variance in the selected games for later suspension days.
\begin{figure}[ht]
    \centering
    \includegraphics[clip,width=0.65\textwidth]{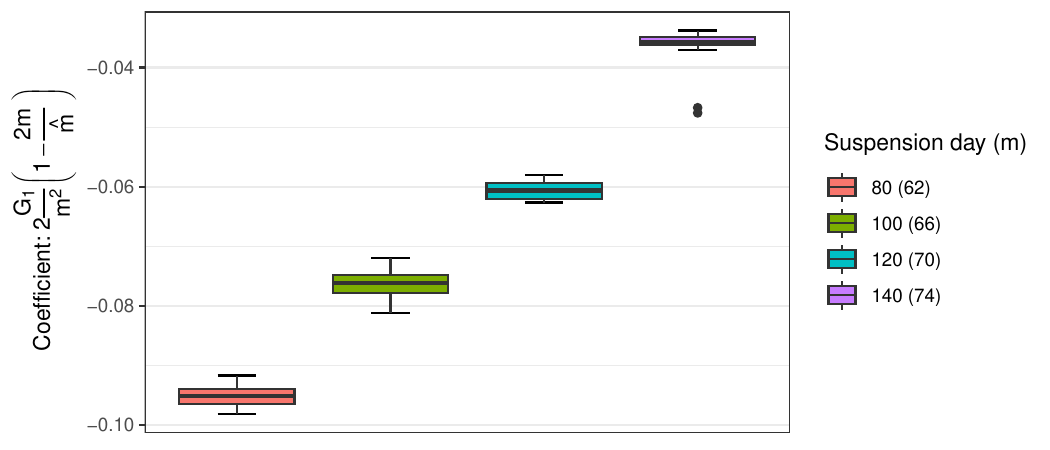}
    \caption{Distribution of the variance coefficients in model PW-DQIP across 14 NBA seasons}
    \label{fig:pd-dqip-var-coeffs}
\end{figure}

Incidentally, we can find a signature of this behaviour aligned with the summarized conclusions in Table \ref{table:short-long-shortened-season} in our computational experiments.
Figure \ref{fig:sharpness-probs-selected-excluded-games-boxplot} illustrates the variance (top row) and the sharpness (bottom row) of the predicted probabilities in i) all remaining games (blue), ii) selected games (green), and iii) excluded games (pink) in our shortened seasons.
As we can see, the sharpness of the selected games for suspension days 80 and 100 is typically lower (and the average variance $\overline{v}$ higher) than that of the excluded games.
Moreover, for later suspension days (i.e., 120 and 140), the pattern for sharpness and variance dissipates due to the fact that less weight is being put on the variance term later in the season (recall Figure \ref{fig:pd-dqip-var-coeffs}).
\begin{figure}[ht]
    \centering
    \includegraphics[clip,width=1\textwidth]{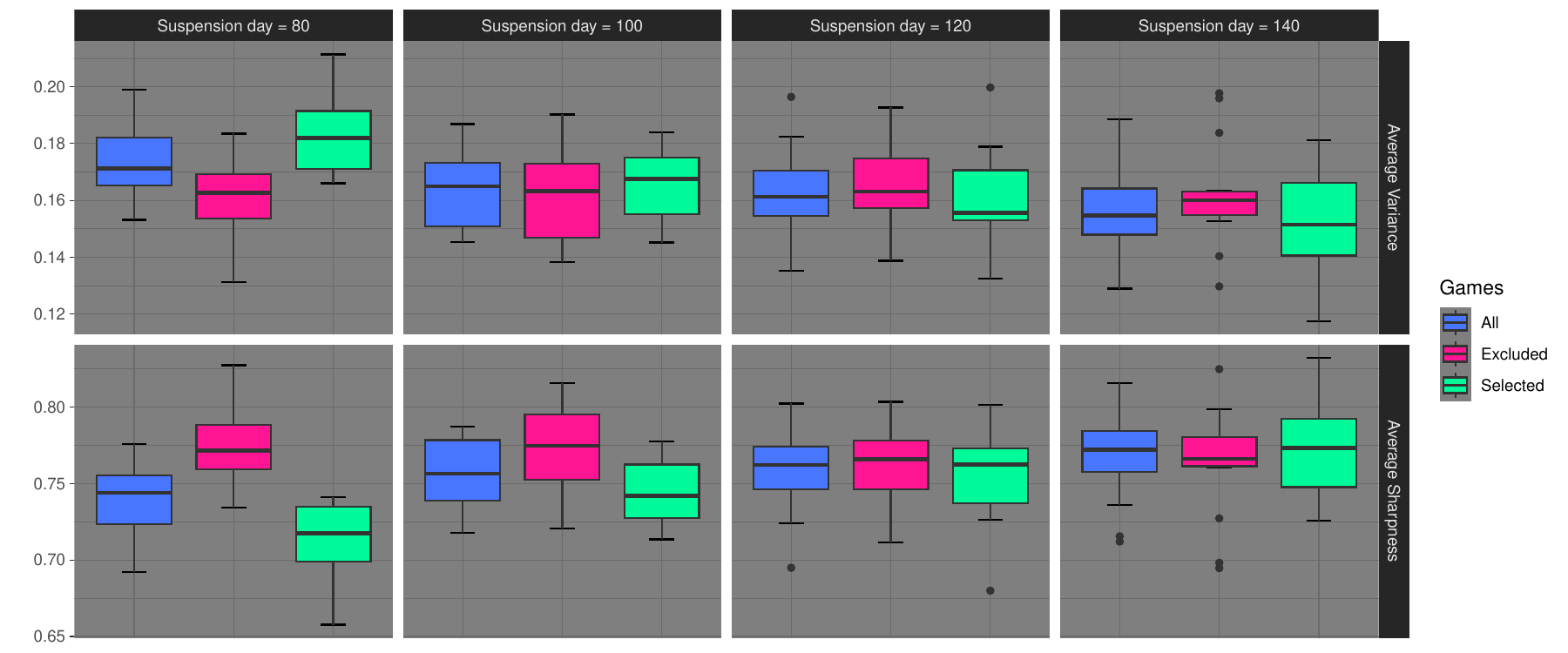}
    \caption{Comparing the average variance and sharpness of the predicted probabilities between selected/excluded games in the shortened season.}
    \label{fig:sharpness-probs-selected-excluded-games-boxplot}
\end{figure}

\section{Frank-Wolfe Algorithm}\label{app:frank-wolfe-algorithm}


Algorithm~\ref{pseudo-code-frank-wolfe} presents the FW algorithm for solving instances of continuous relaxation of PW-DQIP, in which $\bar{X}$ is the continuous relaxation of $X$, and the objective function $f$ is of the form 
\begin{align}
    f(\xx)=\sum\nolimits_{i\in T}\left(\left(\mu_i(\xx)-\hat{\mu}_i\right)^2 +v_i(\xx)\left(1-\frac{2m}{\hat{m}}\right)+ \hat{v}_i \right),
\end{align}
where $\mu_i(\xx)$ and $v_i(\xx)$ are defined in equations \eqref{opt:win-percentage-diff-det-mu} and \eqref{opt:win-percentage-diff-det-v}, respectively. We remark that, since $f$ is a convex quadratic function, its gradient can be computed easily, and a closed-form optimal solution to the line search problem \eqref{fw-ls-app} can be found using the first-order optimality conditions.

\begin{algorithm}[ht]
\caption{PW-FW: Frank-Wolfe algorithm for solving continuous relaxation of PW-DQIP}
\label{pseudo-code-frank-wolfe}
	\begin{algorithmic}[1]
	    \State Let $t\leftarrow 0$, and find an integer solution $\xx^{(0)}\in X$.
	    \State Let $\hat{\xx}^*\leftarrow \xx^{(0)}$
		\While{not converged}
		    \State Compute gradient $d^{(t)}_g=\nabla f_{x_{g}}(\xx^{(t)})$ for each $g\in G$
		    \State Find the integer solution $\hat{\xx}^{(t)}$ by solving the following transportation problem
		    \begin{align}
                \text{[Transportation Problem]}\quad \hat{\xx}^{(t)}=\argmin\nolimits_{\xx\in \bar{X}} \quad & \sum\nolimits_{g\in G}d^{(t)}_g x_g.\label{fw-sub-app}
            \end{align}
            \If{$f(\hat{\xx}^{(t)})<f(\hat{\xx}^*)$}
            \State $\hat{\xx}^*\leftarrow \hat{\xx}^{(t)}$
            \EndIf
            \State Compute the step-size $\gamma^{(t)}$ using the following line search 
            \begin{align}
                \text{[Line Search]}\quad\gamma^{(t)}=\argmin\nolimits_{\gamma \in [0,1]} \quad & f\left((1-\gamma)\xx^{(t)}+\gamma \hat{\xx}^{(t)}\right).\label{fw-ls-app}
            \end{align}
		    \State Update $\xx^{(t+1)}=(1-\gamma^{(t)})\xx^{(t)}+\gamma^{(t)} \hat{\xx}^{(t)}$, and set $t\leftarrow t+1$.
		\EndWhile
	\end{algorithmic}
\end{algorithm}

Given that the atomic solution $\hat{\xx}^{(t)}$ produced by solving the transportation problem \eqref{fw-sub-app} is a feasible integer solution, it provides an upper bound on the optimal value of PW-DQIP. As the algorithm iterates, ${\xx}^{(t)}$ converges to the optimal fractional solution, and $\hat{\xx}^{(t)}$ yields a tighter upper bound. Consequently, the best atomic solution (i.e., $\hat{\xx}^*$) can be used as a near optimal integer solution to PW-DQIP. Henceforth, we refer to this procedure of producing the shortened season $\hat{\xx}^*$ using FW (Algorithm \ref{pseudo-code-frank-wolfe}) as PW-FW.

We would like to remark that our implementation of FW algorithm (Algorithm 1) solves the continuous relaxation of PW-DQIP to optimality, and produces an integer solution as a byproduct thanks to the feasible region of PW-DQIP being totally unimodular.

\textbf{Sub-optimality bounds for Algorithm 1.}
We first note the following property of iterates of FW based on convexity of $f$:
\begin{align*}
    f(\bar{\xx}^*)\ge f(\xx^{(t)}) - \nabla_{\xx}f(\xx^{(t)})^{\top}(\xx^{(t)}-\hat{\xx}^{(t)}),
\end{align*}
which is tight for $\xx^{(t)}=\bar{\xx}^*$. Consequently, we can construct a lower bound on the optimal value of the continuous relaxation of PW-DQIP as
\begin{align*}
    \underline{f} = \max_{t}\left\{f(\xx^{(t)}) - \nabla_{\xx}f(\xx^{(t)})^{\top}(\xx^{(t)}-\hat{\xx}^{(t)})\right\}.
\end{align*}
Let $\hat{\xx}^*$ be the best integer solution produced by FW (i.e., $\hat{\xx}^* = \arg\min_{\hat{\xx}^{(t)}}f(\hat{\xx}^{(t)})$), and $\xx^*$ be the (unknown) integer optimal solution to PW-DQIP. Noting that $\underline{f} \le f(\bar{\xx}^*) \le f(\xx^*) \le f(\hat{\xx}^*)$ we derive the following relative sub-optimality gap 
\begin{align}
    \text{Gap}_{\text{FW}} = \frac{f(\hat{\xx}^*)-\underline{f}}{\underline{f}}, \label{eq:fw-gap-rel}
\end{align}
which upper-bounds the optimality gap with respect to $f(\xx^*)$. More precisely, rewriting the numerator in Eq. \eqref{eq:fw-gap-rel} as
\begin{align}
    f(\hat{\xx}^*)-\underline{f} = \overbrace{f(\hat{\xx}^*)-f(\xx^*)}^{\text{Primal gap}}+ \overbrace{f(\xx^*) - \underline{f}}^{\text{Dual gap}}, \label{eq:fw-gap-abs}
\end{align}
the gap computed in Eq. \eqref{eq:fw-gap-rel} captures the compound effect of quality of the solution produced by FW and  formulation strength of PW-DQIP. Nonetheless, as illustrated in Figure~\ref{fig:FW_gaps}, the sub-optimality gap in Eq. \eqref{eq:fw-gap-rel} is typically very low. For comparison, we also present the absolute gaps as in Eq. \eqref{eq:fw-gap-abs}. 
We remark that our implementation of FW converges to the optimal continuous solution after a few iterations (i.e., $\underline{f}=f(\bar{\xx}^*)$).

\begin{figure}[htbp]
    \centering
    \subfigure[Relative gap]{
    \includegraphics[clip,width=0.35\textwidth]{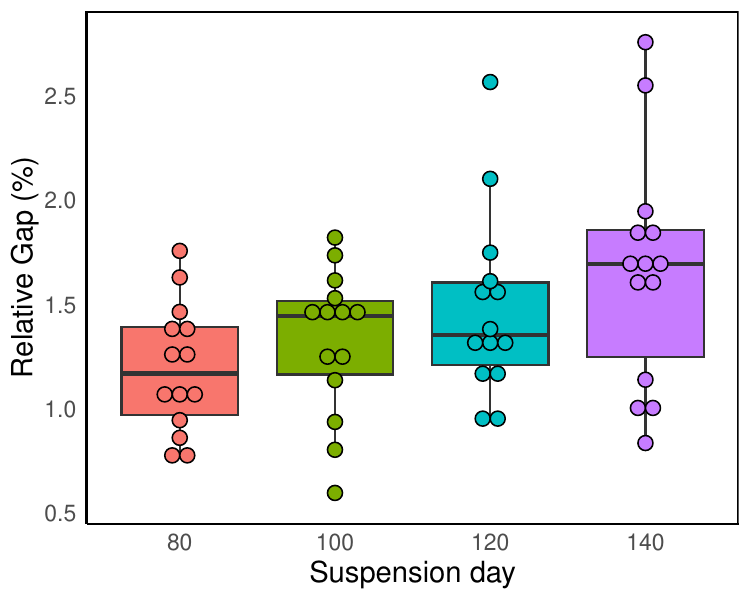}\label{fig:fw_gap_rel}}
    \hspace{0.05\textwidth}
    \subfigure[Absolute gap]{
    \includegraphics[clip,width=0.35\textwidth]{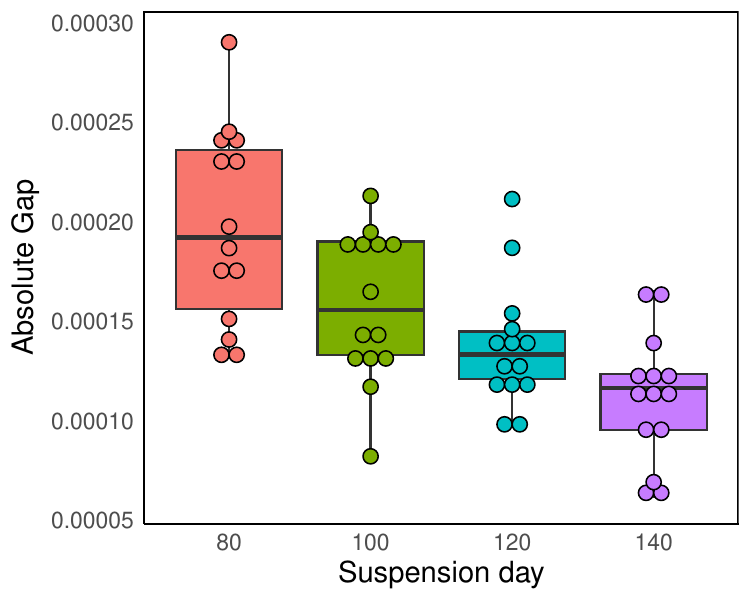}\label{fig:fw_gap_abs}}
    \caption{Upper bounds on the optimality gap of solutions produced by FW. Left: Relative sub-optimality gap (Eq. \ref{eq:fw-gap-rel}). Right: Absolute sub-optimality gap (Eq. \ref{eq:fw-gap-abs}).}
    \label{fig:FW_gaps}
\end{figure}

Finally, we remark that the FW algorithm may be extended to produce near-optimal solutions for PW-MMR as well. However, given the non-smooth min-max (i.e., $\ell_{\infty}$) objective in PW-MMR, a direct implementation of FW may not converge to an optimal continuous solution of PW-MMR \citep[c.f.][for a counterexample]{nesterov2018complexity}; instead, we may minimize a smooth $\ell_p$-approximation of $\ell_{\infty}$ for some large $p$ or minimize the Moreau envelope of the objective function (c.f., \citealp{parikh2014proximal}).

\section{Explanatory Features for Predictive Modeling}\label{app:exp-features}
In this section, explanatory variables are categorized into four groups: overall team performance, basic team-level statistics, advanced team-level statistics, and player-level statistics.

\subsection{Overall Team Performance}\label{app:overall-performance-features}
The most important variable which carries the largest explanatory weight among all the features is \emph{win percentage} of home and guest teams which indicates the relative performance of both teams at the time they play against each other. We also include two variables of the same nature, percentage of home games won by the home team, and percentage of away games won by the guest team, to capture performance variability due to home/away condition. Table \ref{table:performance-features} lists all four overall performance features.
\begin{table}[ht]
\renewcommand{\arraystretch}{1.2}
\centering
\footnotesize
\begin{tabular}{c l}
\toprule
\textbf{Overall Team Performance}    & \multicolumn{1}{c}{\textbf{Definition}}\\ \midrule
\texttt{WPCT} & home team win percentage  \\ 
\texttt{opptWPCT} & guest team win percentage  \\ 
\texttt{WPCTh} & home team win percentage at home  \\ 
\texttt{opptWPCTg} & guest team win percentage on the road \\\bottomrule
\end{tabular}
\vspace{0.2cm}
\caption{Overall Performance Features for home and guest teams.}\label{table:performance-features}
\end{table}

\subsection{Basic Team--Level Statistics}\label{app:team-stats-raw-data}
An obvious choice for an explanatory variable to predict the outcome of basketball games is team-level raw statistics. Over a stretch of games, we can consider average team-level statistics by each team (e.g., average number of points, rebounds, assists, blocks, steals) as explanatory features. Table \ref{table:team-stats-raw} shows the list of such variables for the home team. Note that using a prefix \texttt{oppt} before each variable in Table \ref{table:team-stats-raw} results in the same variable for the guest team, and using a prefix \texttt{diff} for the same set of variables results in the difference between performance of home and guest teams with respect to each variable.
\begin{table}[ht]
\renewcommand{\arraystretch}{1.2}
\centering
\footnotesize
\begin{tabular}{c l}
\toprule
\textbf{Basic Team Features}    & \multicolumn{1}{c}{\textbf{Definition}}\\ \midrule
\texttt{PTS} & Average number of points per game scored  \\ 
\texttt{REB} & Average number of rebounds per game  \\ 
\texttt{AST} & Average number of assists per game \\ 
\texttt{OREB} & Average number of offensive rebounds per game \\ 
\texttt{DREB} & Average number of defensive rebounds per game \\
\texttt{STL} & Average number of steals per game \\ 
\texttt{BLK} & Average number of blocks per game \\
\texttt{TOV} & Average number of turnovers per game \\
\texttt{PF} & Average number of personal fouls per game \\
\texttt{FGM} & Average number of field goals made per game \\
\texttt{FG\%} & Average field goal made percentage \\
\texttt{3PM} & Average number of 3-point field goals made per game \\
\texttt{3P\%} & Average 3-point shot percentage \\
\texttt{FTM} & Average number of free throws made per game \\
\texttt{FT\%} & Average free throw made percentage \\
\texttt{PITP} & Average number of points in the painted area per game \\
\texttt{FBPs} & Average number of fast--break points per game \\
\texttt{2ndPTS} & Average number of second chance points per game \\
\texttt{PTSOFFTO} & Average number of points off of opponent's turnovers per game \\
\texttt{Poss} & Average number of possessions per game \\\bottomrule
\end{tabular}
\vspace{0.2cm}
\caption{Team--level statistics used as explanatory features.}\label{table:team-stats-raw}
\end{table}

\subsection{Advanced Team--Level Statistics}\label{app:team-stats-advanced-data}
These are advanced features calculated based on the raw data shown in Table \ref{table:team-stats-raw}. The goal of introducing and using these advanced features is to highlight strengths and weaknesses of each team adjusted by their style of play (e.g., reliance of each team on 3-point shots, defensive style of play). For instance, \texttt{FG\%} and \texttt{3P\%} are two basic statistics while \emph{effective field goal percentage}, denoted by \texttt{eFG\%}, computes a weighted average field goal percentage, applying a weight of 2 to regular field goals and a weight of 1 to 3-point shots, scaled by the number of field goal attempts.
Table \ref{table:team-stats-advanced} contains the list of advanced team-level statistics. The same set of variables are defined for the guest team (identified by prefix \texttt{oppt}) and the difference of each variable between the two teams (identified by the prefix \texttt{diff}).
\begin{table}[ht]
\centering
\footnotesize
\begin{tabular}{c p{12cm}}
\toprule
\textbf{Advanced Team Features}    & \multicolumn{1}{c}{\textbf{Definition}}\\ \midrule
\texttt{OffRtg} & Offensive rating, which is the number of points scored per 100 possessions  \\ 
\texttt{DefRtg} & Defensive rating, which is the number of points allowed per 100 possessions  \\ 
\texttt{NetRtg} & Net rating of a team is calculated by subtracting \texttt{DefRtg} from \texttt{OffRtg} \\ 
\texttt{AST\%} & An estimate of the percentage of field goals assisted by team players per game \\ 
\texttt{AST/TO} & Assists per turnover ration, which is the number of assists per team divided by the number of turnovers the team has committed in a game \\
\texttt{ASTRatio} & Average number of assists per 100 possessions \\ 
\texttt{OREB\%} & Offensive rebound percentage, which is an estimate of the percentage of available offensive rebounds a team grabs per game \\
\texttt{DREB\%} & Defensive rebound percentage, which is an estimate of the percentage of available defensive rebounds a team grabs per game \\
\texttt{REB\%} & Total rebound percentage which is an estimate of the percentage of total available rebounds a team grabs per game \\
\texttt{TOV\%} & Turnover percentage, which is the percentage of plays that end in a player or team turnover \\
\texttt{eFG\%} & Effective field goal percentage, which measures field goal percentage adjusting for made 3-point field goals being 1.5 times more valuable than made 2-point field goals. \\
\texttt{TS\%} & True shooting percentage, a measure of shooting efficiency which differentiates between the number of points awarded by a regular field goal, a 3-point field goal, and a free throw. \\\bottomrule
\end{tabular}
\vspace{0.2cm}
\caption{Advanced Team--level statistics calculated based on raw features from Table \ref{table:team-stats-raw}.}\label{table:team-stats-advanced}
\end{table}

\subsection{Player-Level Statistics}\label{app:player-stats}
Each of the features introduced in \cref{app:team-stats-raw-data} can be defined for an individual player as well. With some adjustments, all the features introduced in \cref{app:team-stats-advanced-data} can also be defined for individual players. Given there are 15 players on the roster for any NBA team, with at least 9 playing considerable minutes each night, the total number of features will grow substantially large, should we choose to define player-level features corresponding to team-level features in Tables \ref{table:team-stats-raw} and \ref{table:team-stats-advanced}. To tackle this issue, there are alternative ways to represent the efficiency of individual players using a combination of raw data. \emph{Efficiency} rating introduced by \cite{efficiency-rating} is one way to combine individual statistics into a single number. The formula is the following:
\begin{equation}
    \text{\texttt{EFF}} = \text{\texttt{PTS}} + \text{\texttt{REB}} + \text{\texttt{AST}} + \text{\texttt{STL}} + \text{\texttt{BLK}} - \text{Missed \texttt{FG}} - \text{Missed \texttt{FT}} - \text{\texttt{TOV}}\label{eq:eff-formula}
\end{equation}

We calculate the \texttt{EFF} rating for each player according to \eqref{eq:eff-formula} (using all the games prior to the game under study), and we represent the home and guest teams by their top 10 players, sorted according to their \texttt{EFF} values. Let $\text{\texttt{hp}}_i$ denote the \texttt{EFF} rating for the $i^{\text{th}}$ best player of the home team, and let $\text{\texttt{gp}}_i$ denote the \texttt{EFF} rating for the $i^{\text{th}}$ best player of the guest team, $i \in \{1, 2, \dots, 10\}$. We can use the mean and standard deviation of these 10 numbers to represent overall efficiency of players on each roster and the discrepancy of \texttt{EFF} ratings between players. Let (\texttt{EFF-mean}, \texttt{EFF-std}) and (\texttt{opptEFF-mean}, \texttt{opptEFF-std}) represent the average \texttt{EFF} and standard deviation of \texttt{EFF} values for home and guest teams, respectively.

Lastly, we would like to highlight the importance of including player-level statisics in our model's application to suspension scenarios such as the NBA lockouts, where no games were played initially. To address this, we propose a modified approach using player-level statistics from previous seasons to create team-level features. This method accounts for off-season player movements, ensuring our predictions remain relevant for the upcoming season. By aggregating individual player performances, such as points and assists, we can effectively simulate team capabilities despite the absence of current season games. This adaptation demonstrates our model's flexibility and its ability to provide accurate predictions in a variety of scenarios, including those without any pre-season games.

~

\section{Supplementary Results}

\subsection{Cross Simulation}\label{app:cross_sim}
Figure~\ref{fig:cross_simulation} illustrates the simulation results for different choices of predictive models in the prescriptive and simulation phases. Each point corresponds to the concordance value between the shortened and full seasons averaged over 14 seasons.
\begin{figure}[ht]
    \centering
    \includegraphics[width=\textwidth]{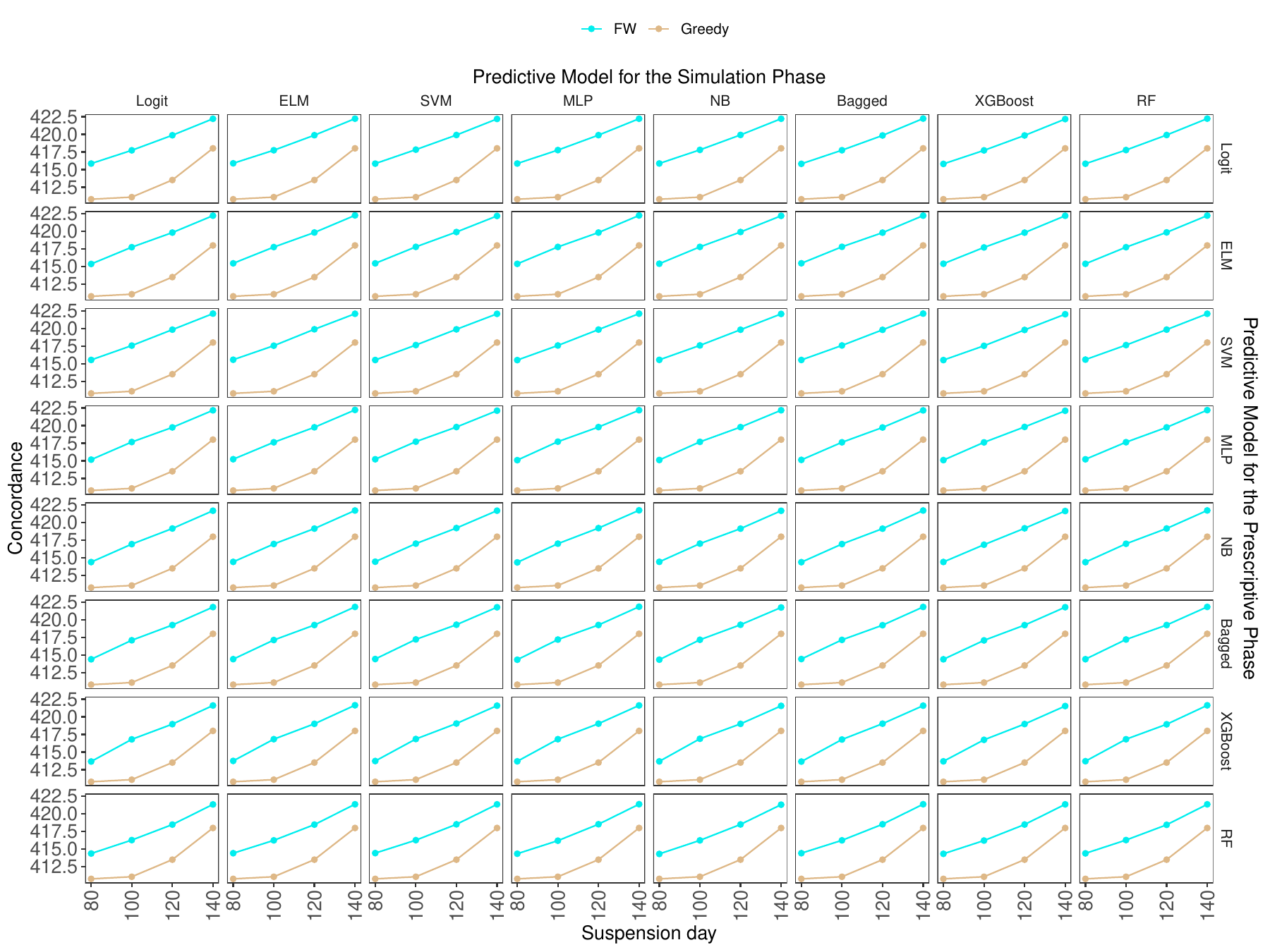}
    \caption{Cross comparison of simulation results based on different choices of predictive models used for the prescriptive phase (rows) and predictive models used for the simulation phase (columns).}
    \label{fig:cross_simulation}
\end{figure}

\subsection{Practical Implications: Comparing PW-FW, Greedy the Status Quo Models}\label{app:practical-implications-boxplots}
In \cref{sec:presc-model-practical-impl}, we introduced three metrics to gauge the practical implications of various shortened season plans and we compared our best performing prescriptive model (i.e., PW-MMR) with two benchmarks: Greedy and Status Quo in overall agreement. In this section, we plot the agreement distributions for each of the 4 suspension dates for models PW-MMR and the baseline Greedy in Figure \ref{fig:success_rates_greedy_boxplots}, and similar distributions for PW-MMR and the baseline Status Quo in Figure \ref{fig:success_rates_Status_Quo_boxplots}.
Our best-performing solution method, PW-MMR, outperforms both baseline solutions in all three success rate metrics both in terms of average percentage and the variation across 14 NBA seasons.
Similar to the simulation results presented in Figure \ref{fig:sim_results}, as the suspension day increases, the margin of improvement with respect to the Status Quo model gets smaller.

\begin{figure}[ht]
    \centering
    \includegraphics[clip,width=\textwidth]{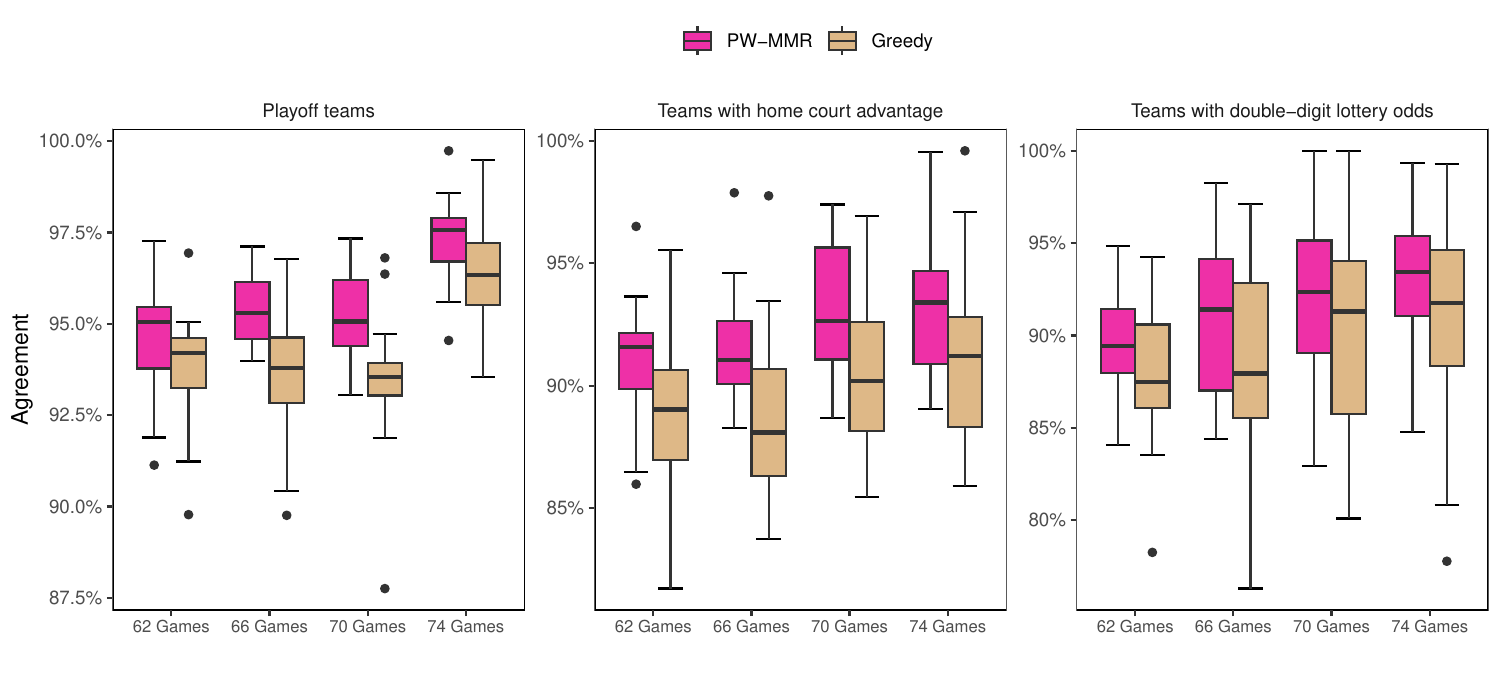}
    \caption{Agreement distributions for PW-MMR and Greedy.}
    \label{fig:success_rates_greedy_boxplots}
\end{figure}

\begin{figure}[ht]
    \centering
	\includegraphics[clip,width=\textwidth]{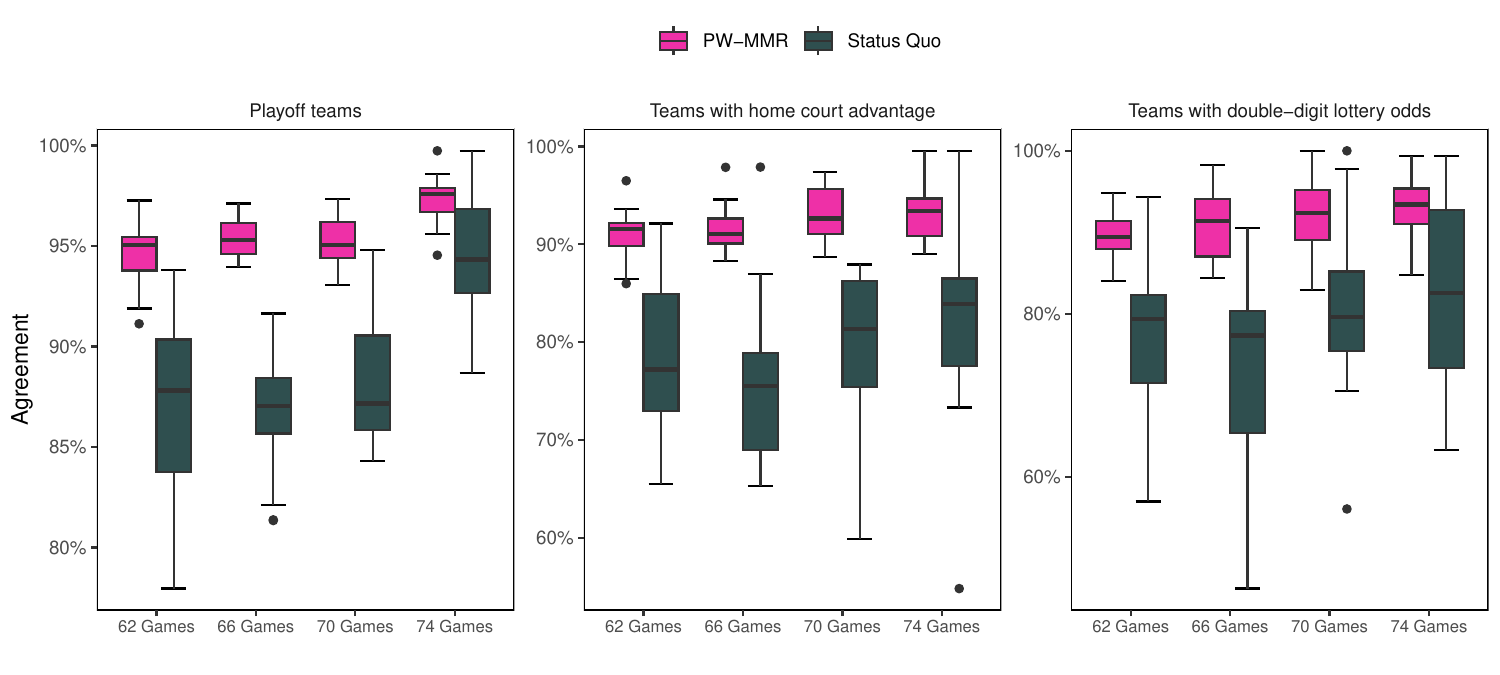}
    \caption{Agreement distributions for PW and Status Quo.}
    \label{fig:success_rates_Status_Quo_boxplots}
\end{figure}

\subsection{Strength of Schedule Extension: Comparison Across Seasons}\label{app:SoS-lineplots}
In \cref{sec:sos}, we introduced model PW-SoS which amends the model PW-DQIP (presented in \cref{sec:prescriptive}) by adding a constraint making sure that the strength of schedule discrepancy (SSD) in the shortened season is no larger than $\epsilon \%$ compared to the SSD metric in the full season. In our experiments, we tested 4 different choices of the parameter $\epsilon$ including (0.02, 0.03, 0.05, 0.1). Figure \ref{fig:OWP_overall_ssd_sim} in \cref{sec:sos} illustrates the comparison between PW-SoS models as well as PW and Greedy using values averaged across 14 seasons. The following two figures in this section present the same comparison in terms of the metric SSD for individual seasons as lineplots. Figure \ref{fig:OWP_difference_PWS_PW_Greedy} plots SSD according to PW and the our benchmark Greedy as well as the best choice of PW-SoS (with $\epsilon = 0.02$). Figure \ref{fig:OWP_difference_PWS_levels} plots the SSD values in individual seasons for 4 choices of $\epsilon$. We can easily conclude that the model PW-SoS with $\epsilon = 0.02$ has the lowest overall SSD, thus performing better than the other three.

\begin{figure}[t!]
    \centering
    \includegraphics[clip,width=\textwidth]{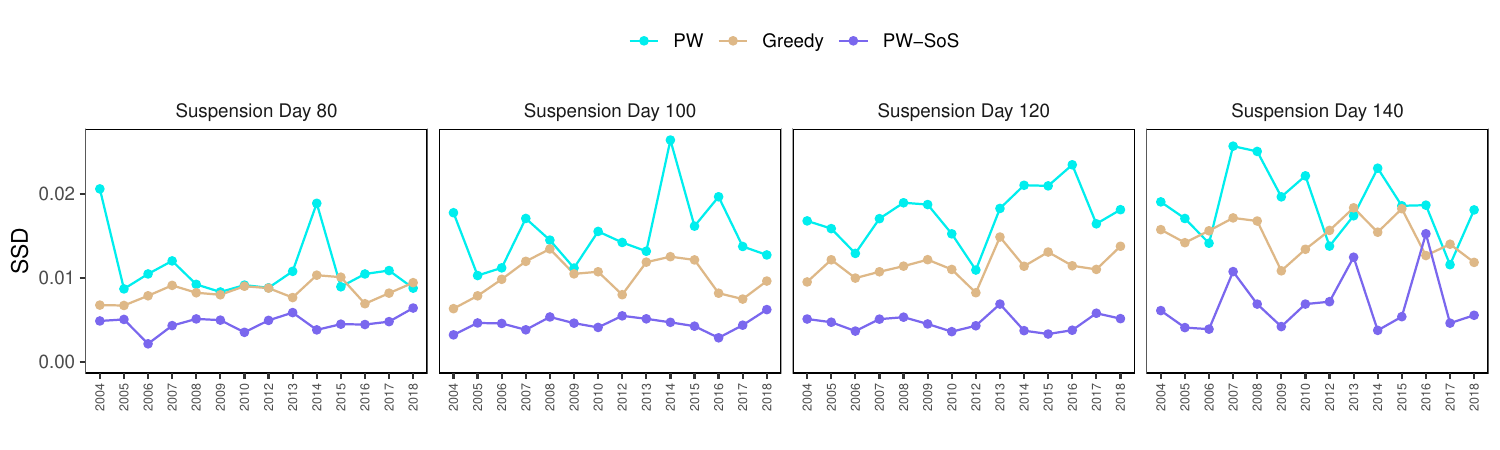}
    \caption{Strength of Schedule Discrepancy (SSD) between shortened season and full season solutions based on their post-suspension strength of schedule, comparing the distribution of SSD over 14 NBA seasons in PW-SoS, PW and Greedy.}
    \label{fig:OWP_difference_PWS_PW_Greedy}
\end{figure}

\begin{figure}[!ht]
    \centering
    \includegraphics[clip,width=\textwidth]{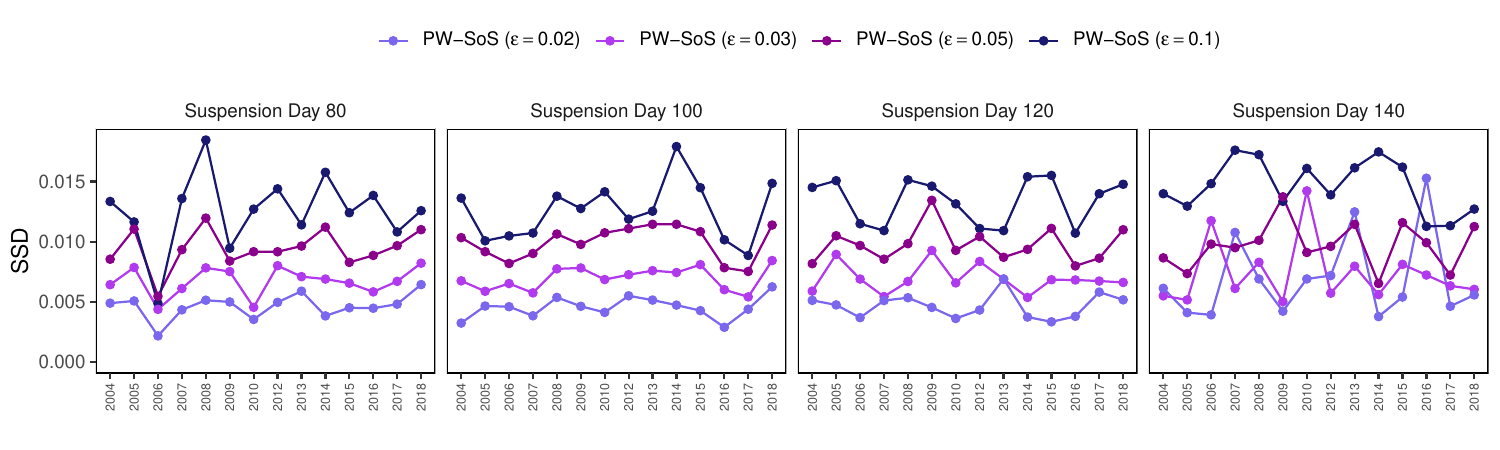}
    \caption{Strength of Schedule Discrepancy (SSD) between shortened and full season solutions based on their post-suspension strength of schedule, comparing the distribution of SSD over 14 NBA seasons for PW-SoS with various choices for parameter $\epsilon$.}
    \label{fig:OWP_difference_PWS_levels}
\end{figure}

Note that we tested the results according to a few other definitions of strength of schedule, including Relative Percentage Index (RPI) as defined in \cite{NBAStufferSoS}, and the same conclusion stands.

\subsection{Suggestions for the 2019--20 Season}\label{app:suggestions}

In this section, we present the results of our two-phase analytics approach applied to the 2019--20 NBA regular season which was suspended on March 11, 2020 due to the COVID-19 pandemic. We consider 74 games per team as the target length of the shortened season, thus canceling 8 games per team from the remainder of the season. As a result, out of 259 remaining games, we select 139 games to be played in the shortened season using the PW-FW model. The set of selected and canceled games are shown in Figure \ref{fig:calendar74_homeaway}.

\begin{figure}[t!]
    \centering
	\includegraphics[clip,width=1\textwidth]{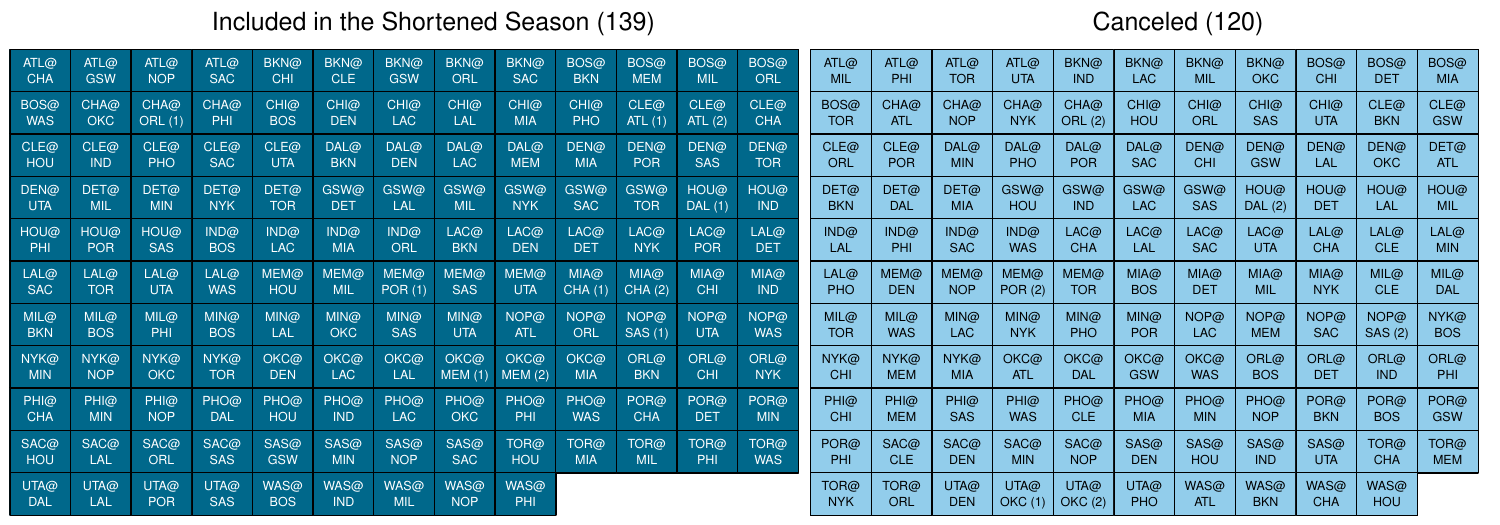}
    \caption{Selected/canceled games for the remainder of the 2019--20 season according to PW-FW.}
    \label{fig:calendar74_homeaway}
\end{figure}

According to the NBA's resumption plan announced on June 26, 2020~\citep{nba-return}, 22 teams (the top 13 teams from the west and top 9 teams from the east) were invited to Orlando, Florida to play 8 more games each to conclude the 2019--20 NBA regular season. In effect, the invited teams will have a shortened season ranging from 71 to 75 games in total, and this variation is a consequence of some teams having played a few more games than others as of the suspension date.
We compare the shortened season plan shown in Figure \ref{fig:calendar74_homeaway} to the NBA's resumption plan. Using the same Monte Carlo simulation approach described in \cref{sec:presc-model-results}, we evaluated our proposed solution and the NBA's return plan using 1,000 scenarios and computed the average concordance relative to the full season ranking. The concordance for our proposed solution is 404.7, while it is 392.25 for the NBA's resumption plan. As a result, on average, our proposed solution predicts the relative positioning of at least 12 additional pairs of teams correctly, compared the the NBA's return plan.



\end{APPENDICES}

\end{document}